\documentclass[11pt]{article}
\usepackage{amsfonts,amssymb, graphicx,a4,bm,bbm}
\usepackage{float}
\usepackage{soul}
\usepackage[normalem]{ulem}
\usepackage{color}

\newcommand{\zeq}{\setcounter{equation}{0}}

\newtheorem{theo}{Theorem}[section]
\newtheorem{lem}[theo]{Lemma}

\newtheorem{defi}[theo]{Definition}

\definecolor{Red}{cmyk}{0,1,1,0}

\let\a=\alpha \let\b=\beta \let\ch=\chi  
 \let\g=\gamma     \let\k=\kappa 
\let\m=\mu  \let\o=\omega    \let\p=\pi \let\ph=\varphi
\let\r=\rho \let\s=\sigma \let\t=\tau 
 \let\x=\xi 
\let\D=\Delta \let\F=\Phi \let\G=\Gamma \let\L=\Lambda 
\let\O=\Omega    

\let\\=\noindent
\def\ee{\end{equation}}
\def\be{\begin{equation}}
\def\vv{\vskip.2cm}
\let\a=\alpha \let\b=\beta \let\ch=\chi  
 \let\g=\gamma     \let\k=\kappa 
\let\m=\mu  \let\o=\omega    \let\p=\pi \let\ph=\varphi
\let\r=\rho \let\s=\sigma \let\t=\tau 
 \let\x=\xi 
\let\D=\Delta \let\F=\Phi \let\G=\Gamma \let\L=\Lambda 
\let\O=\Omega    

\let\\=\noindent
\def\ee{\end{equation}}
\def\be{\begin{equation}}
\def\vv{\vskip.2cm}

\def\rfp{R^{\mbox{\tiny\rm FP}}}

\def\phd{\ph^{\mbox{\tiny\rm D}}}

\def\phkp{\ph^{\mbox{\tiny\rm KP}}}
\def\VU{{\mathbb{V}}}\def\0{\emptyset}\def\Z{\mathbb{Z}}
\def\fa{{\mathfrak{F}}}
\def\obj{{\rm supp}}
\def\prob{{\rm Prob}}
\def\ev{{\mathfrak{e}}}

\def\PP{{\mathcal P}}

\def\FF{{\mathfrak{F}}}

\def\GG{{\cal G}}

\def\0{\emptyset}
\def\Om{\Omega}
\title{Moser-Tardos  resampling algorithm, entropy compression method and the subset gas}
\author{
\\
Paula M. S. Fialho, Bernardo N. B. de Lima, Aldo Procacci\\
\\
\small{ Departamento de Matem\'atica UFMG}
\small{ 30161-970 - Belo Horizonte - MG
Brazil}}

\begin{document}

\maketitle

\begin{abstract}
We establish a connection between  the entropy compression method and  the Moser-Tardos  algorithmic version
of the Lov\'asz local lemma through  the cluster expansion of the  subset gas. We also show that the Moser-Tardos resampling algorithm and the entropy compression bactracking algorithm  produce identical bounds.
\end{abstract}

\vskip.2cm
{\footnotesize

}
\vskip.5cm

\vskip.5cm

\section{Introduction}

\subsection{The Lov\'asz Local Lemma}
The Lov\'asz Local Lemma (LLL), originally formulated by Erd\"os and Lov\'asz in  \cite{EL}, is a powerful tool in the framework of the
probabilistic method  used in an  impressive quantity of applications in combinatorics such as graph coloring problems,
K-sat, latin transversal, etc.,   (see \cite{AS} and references therein for a review). Its basic idea is to prove  the  existence of some combinatorial object with certain desired property
(e.g. such a proper coloring of  the vertices of a graph) by identifying  a family $\FF$  of (bad) events in some probability space $\O$
whose presence, even of only one of them, spoils  the object under analysis and whose
simultaneous non-occurrence  guarantees that the object under analysis is actually present.
Denoting by  $\overline \ev$  the complement
event of $\ev\in \fa$,  the Lov\'asz local lemma provides a condition on  the probabilities $\prob(\ev)$  in order to ensure that
$\prob(\bigcap_{\ev\in \fa}\overline \ev)>0$.  To formulate explicitly this condition
one need to identify a so-called dependency graph for the family $\fa$. That is to say,  a graph $\GG$ with vertex set $\FF$ and edge set
such that each event $\ev\in \fa$   is independent of
the $\s$-algebra generated by the collection of events $\fa \setminus \Gamma^*_\GG(\ev)$, where $\Gamma_\GG^*(\ev)=\G_\GG(\ev)\cup\{\ev\}$ and  $\G_\GG(\ev)$
is the set of all the events of $\fa$   adjacent to $\ev$ in $\GG$. According to the usual terminology, $\Gamma^*_\GG(\ev)$ is called the neighborhood of $\ev$ in $\GG$ while  $\G_\GG(\ev)$ is called the {\it punctured neighborhood} of $\ev$ in $\GG$.

\\Once the dependency graph
of the family $\fa$  has been determined, the LLL can be stated as follows.

\begin{theo}[Lov\'asz local Lemma]\label{lov}  Let $\fa$ be a finite family of events in a probability space  $\Om$ and let  $\GG$ be a dependency graph for  $\FF$.  Let  $\bm{\mu}=\{\mu_\ev\}_{\ev\in \FF}$ be a collection of non-negative numbers. If, for each $\ev\in \fa$,
\begin{equation}
\label{cond}
\prob(\ev)\;\leq\;  {\mu_\ev \over \prod_{\ev'\in \G^*_\GG(\ev)}(1+\m_{\ev'})},
\end{equation}
then
\be\label{thesis}
\prob\Big(\bigcap_{\ev\in \FF}\overline \ev\Big)\,> \,0.
\ee
\end{theo}

\subsection{The abstract polymer gas}

\\The abstract polymer system (APS) is a discrete model originally proposed by Kotecky and Preiss \cite{KP}
as a generalization of a lattice polymer  model introduced by Gruber and Kunz \cite{GK} in 1968. Its relevance in statistical mechanics is very important since it is a widely used  tool to analise a large number of systems in physics, such as discrete spin systems, continuous and discrete
particle systems, percolative models and even quantum field theories.

\\The APS is  defined by a triple $(\PP,\mathbf{w}, W)$ where $\PP$ is a countable (possibly infinite) set whose elements are called polymers, $\mathbf{w} :\PP\to \mathbb{C}$ is a function which associates to each polymer  $\g\in \PP$ a  complex number $w_\g$,
called the activity of the polymer $\g$, and  $W:\PP\times \PP\to \{0,1\}$ is a  function called the {\it Boltzmann factor}, such that $W(\g,\g)=0$ and $W(\g,\g')=W(\g',\g)$ for all $\{\g,\g'\}\subset \PP$.
Usually the pair $\{\g,\g'\}$ is  called {\it incompatible} when $W(\g,\g')=0$ and {\it compatible} when $W(\g,\g')=1$.

\\Let $\GG$ be the simple graph with vertex set $\PP$ and edge set  formed  by the pairs $\{\g,\g'\}\subset \PP$ such that
$W(\g,\g')=0$. The graph $\GG$, which is uniquely determined by the Boltzmann factor
$W$, is sometimes called the {\it support graph} of $W$.
The  neighborhood of the vertex $\g$ in the graph $\GG$ is the set
$\G^*_\GG(\g)=\{\g'\in \PP: ~W(\g,\g')=0\}$ formed by all polymers incompatible with $\g$.
An independent set  of the support graph $\GG$ is    a set $Y$ of polymers  such that each pair $\{\g,\g'\}\subset Y$ is compatible. We denote by $I(\GG)$ the set formed by all finite independent sets of $\GG$.

\\Given a finite collection of polymers $\L\subset \PP$, the grand canonical partition function of the APS  at ``finite volume"  $\L$
is given by
$$
Z_\L(\bm w)= \sum_{S \subset \L \atop S \in I(\GG)}{\prod_{\g \in S}{w_\g}}.
$$
This is a key quantity  since the thermodynamic properties of the system can be derived from it. In particular,
a fundamental question  physicists are interested in,  is
to find radii $\mathbf{R} = \{R_\g\}_{\g \in \PP}$ (with $R_\g\ge 0$ for all $\g\in \PP$) such that  the partition function  $Z_\L(\bm w)$, for any $\L$ finite,
%and the correlations (\ref{corr})
is free of zeros  for all  complex activities $\bm w$ within the polydisk $\{|w_\g|<R_\g\}_{\g\in \PP}$ (shortly $\bm w\le \mathbf{ R}$).
This would guarantee that the logarithm of the partition function, which is related to the pressure of the system, is analytic in such regions, so that no phase transitions occur.
The best current lower bound for such radii $\mathbf{R}$ is   due to Fern\'andez and Procacci \cite{FP} who improved the older bounds
due to Kotecky and Preiss \cite{KP} and Dobrushin \cite{Do}  proving the following theorem.

\begin{theo}[Fern\'andez-Procacci criterion]\label{FPC}
Let $\bm\mu= \{\mu_\g\}_{\g \in \PP}$ be a collection of nonnegative numbers such that
 \be\label{fp}
 |w_\g| \leq \rfp_\g\equiv \frac{\mu_\g}{\Xi_\g(\bm \mu,\GG)},~~~~~~~~~~\forall \g\in \PP
 \ee
 with
\be\label{Xi}
\Xi_\g(\bm \mu,\GG)=~\sum_{S\subseteq\Gamma_\GG^*(\g)\atop S\in I(\GG) } \prod_{\g'\in S}\mu_{\g'}.
\ee
Then, for all finite $\L\subset \PP$, $Z_\L(\mathbf{w})\neq 0$.
\end{theo}
The Kotecky-Preiss and the Dobrushin criteria can be formulated  analogously with the only difference that function $\Xi_\g(\bm \mu,\GG)$ appearing in the r.h.s. of (\ref{fp})
is replaced respectively  by

\be\label{Kope}
{\phkp_\g(\bm\mu)}=%~\sum_{S\subseteq\Gamma^*_\GG(\g)} \prod_{\g'\in S}\mu_{\g'}=
e^{\sum_{\g'\in \G^*_\GG(\g)}\m_{\g'}},
\ee
and
\be\label{Dobro}
{\phd_\g(\bm\mu)}=~\sum_{S\subseteq\Gamma^*_\GG(\g)} \prod_{\g'\in S}\mu_{\g'}=\prod_{\g'\in \G^*_\GG(\g)}(1+\m_{\g'}).
\ee
The bound on radii $\mathbf{R}$ given by the Fern\'andez-Procacci criterion (\ref{fp})
is always greater than the
bounds on the same radii  given by  the Kotecky-Preiss and the Dobrushin criteria since
$$
\exp\Big\{\sum_{\g'\in \G^*_\GG(\g)}\m_{\g'}\Big\}\ge \prod_{\g'\in \Gamma^*_\GG(\g)} (1+\mu_{\g'})= \sum\limits_{S\subseteq  \Gamma^*_\GG(\g)}
\prod_{\g'\in S}\mu_{\g'} \ge
~\sum_{S\subseteq\Gamma_\GG^*(\g)\atop S\in I(\GG) } \prod_{\g'\in S}\mu_{\g'}.
$$
%\\It is worth to mention that the improved Fern\'andez-Procacci criterion for the APS
%has been recently  employed in some applications in statistical  mechanics \cite{FPS}, \cite{MP}.

\subsection{The connection between the LLL and the APS}
It is important  to remark  that the LLL  criterion (\ref{cond}) is a {\it sufficient} condition in order to the thesis (\ref{thesis}) to hold. In 1985, Shearer \cite{Sh} presented  a necessary and sufficient condition for (\ref{thesis}) to hold. The Shearer condition was actually constituted by a set of several conditions
which were very difficult (if not impossible)  to be checked in practical applications. Probably for this reason Shearer's result went somehow  overseen until 2005. In this year Scott and Sokal \cite{SS}, inspired by Shearer's work, showed  that there was a quite surprising connection between the Lov\'asz Local Lemma   and  the  abstract polymer gas.
Scott and Sokal showed that, given the family of events $\FF$ and their dependency graph $\GG$, the Shearer criterion is  equivalent to require that the probabilities  of the bad events,
$\prob(\ev)$, fall in the zero-free  region of the partition function of the APS whose support graph coincides with the dependency graph $\GG$
of the family $\FF$.  So, once rephrased  in the statistical mechanics lingo,   it is no surprise that
the Shearer criterion was unusable  in practice. On the other hand Scott and Sokal observed that this equivalence implies that the  LLL criterion (\ref{cond})
coincides with the aforementioned Dobrushin criterion.  Later,
Bissacot et al. \cite{BFPS}, via  the connection disclosed in \cite{SS} and Theorem \ref{FPC},    improved  the LLL criterion (\ref{cond})  as follows.

\begin{theo}[Cluster expansion local lemma (CELL)]\label{bis}
With the same hypothesis of the Theorem \ref{lov},  if, for each event  $\ev\in \fa$
\begin{equation}
\label{cond1n}
\prob(\ev)\;\leq\;  {\mu_\ev \over  \Xi_\ev(\bm \m,\GG)},
\end{equation}
with
\be\label{Xi2}
\Xi_\ev(\bm \m,\GG)=\sum_{S\subseteq\Gamma_\GG^*(\ev)\atop S\in I(\GG)} \prod_{\ev'\in S}\mu_{\ev'},
\ee
then
$$
{\prob}\Big(\bigcap_{\ev\in \fa}\bar \ev\Big)\,> \,0.
$$
\end{theo}
As observed above this is clearly an improvement w.r.t. Theorem \ref{lov} since
$$
\sum_{S\subseteq\Gamma_\GG^*(\ev)\atop S\in I(\GG)} \prod_{\ev'\in S}\mu_{\ev'}\le \sum_{S\subseteq\Gamma_\GG^*(\ev)} \prod_{\ev'\in S}\mu_{\ev'}= \prod_{\ev'\in \G^*_\GG(\ev)}(1+\m_{\ev'}).
$$
Condition (\ref{cond1n}) has been shown to be effective  in several applications of the LLL (see e.g. \cite{NPS} and \cite{BKP}). This new criterion is nowadays known as
``Cluster Expansion (CE) criterion".
\vskip.2cm

\\Formulas (\ref{cond}) and  ( \ref{cond1n}) on one hand and formulas (\ref{fp}) and (\ref{Dobro}) on the other hand  show in a crystal way   the evident connection between the  LLL and the APS:
\begin{itemize}
\item[-] events $\ev\in \FF$ in LLL correspond to polymers $\g\in \PP$ in APS;

\item[-] dependents events  correspond to incompatible polymers;

\item[-] probabilities of events correspond to   the absolute values of activities of polymers;

\item[-] probability of the good event to be positive corresponds to require that  the partition function evaluated at  $-|\bm w|$ (which is the worst case, see e.g. \cite{SS} or \cite{BFP}) to be strictly positive.
\end{itemize}

\subsection{Moser-Tardos algorithmic version of the LLL}
\\The LLL is very general, in particular in   its   statement and  proof  there is no need to specify anything about the probability space. Of course, in the applications
the probability space $\Om$ has to be specified and it is natural to wonder, once condition (\ref{cond}) is satisfied, if it is possible to find a  polynomial algorithm  in this specified probability space  able to find a configuration in $\Om$ which
realizes the event  $\bigcap_{\ev\in \FF}\overline \ev$.
During many years researchers have tried to find  methods to devise  general  algorithms able to find such a  configuration for  as many as possible applications covered by the LLL. These efforts have been only partially successful in the sense that the class  of example for which an efficient algorithm could be found  was limited  and
 the condition (\ref{cond}) got  worse, see for example \cite{A},\cite{B}. Such situation changed radically in 2009 when in
a  breakthrough paper \cite{MT}  Moser and Tardos presented a fully algorithmic version of the LLL which covered the vast majority of LLL applications.

\\The scheme   proposed  by    Moser and Tardos  is called nowadays the {\it  variable setting}. It starts from the assumption
that  the probability space $\O$  is a product space generated by
a collection of  {\it mutually independent} random variables  $\{\psi_x\}_{x\in \L}$, where  $\L$ is  a finite set whose elements will be called {\it atoms} hereafter.
In general, each random variable $\psi_x$ takes values in its own space $\Psi_x$ according to its own distribution, but in a vast majority of the applications these
 variables take values in a common finite space $\Psi$ (e.g.  the set of colors). Therefore $\O=\prod_{x\in \L}\Psi_x\equiv \Psi_\L$
 and an element $\o\in \Psi_\L$ is called a {\it configuration}.
An event
$\ev$ in such a probability space $\O$  (also called a {\it flaw}) is a subset of $\Psi_\L$ and it depends in general  on all the random variables
$\{\psi_x\}_{x\in \L}$ generating $\O$. When $\ev$ depends only on variables
of a proper subset $A\subset \L$ (i.e. $\ev$ is fully determined by variables $\Psi_A=\{\psi_x\}_{x\in A}$), we say that the event $\ev$ is {\it tempered}
and we write ${\rm supp}(\ev)=A$  and ${\rm supp}(\ev)=A\subset \L$ is sometimes called the
{\it scope} of $\ev$. A tempered event $\ev$ is hereafter called  {\it elementary} if it is formed only by one single configuration $\o\in \Psi_{\obj(\ev)}$.

  \\Moser and Tardos considered situations in which all of the bad events constituting the family $\FF$ are tempered.
With this assumption   any two events $\ev, \ev'$ of $\fa$ such that $\obj(\ev)\cap \obj(\ev')= \emptyset$
are necessarily independent.
This implies that
the graph $\GG$  with vertex set $\fa$ and edge set  constituted by the pairs
$\{\ev,\ev'\}$ such that $\obj(\ev)\cap \obj(\ev')\neq \emptyset$ is  a natural dependency graph for the family  $\fa$.

\vskip.2cm

\\In this setting  Moser and Tardos defined the following algorithm.
%\vskip.2cm
%{\small \\\textsc{Resampling}.
%\vskip.1cm
%\\- Step 0. Sample a random  configuration $\o\in \Psi_\L$ and set $\o_0=\o$ as the configuration at step $0$.
%
%\\- Step $i$ (for $i\ge 1$). Let $\o_{i-1}\in \Psi_\L$ be the configuration at step $i-1$. If $\o_{i-1}$ is such that  some event
%of the family  $\fa$ occurs,
%select one  of them (at random or according to some deterministic rule), say $\ev$,
%   and  choose an  evaluation (resampling) $\o'\in \Psi_{\obj(\ev)}$. Set
%$\o_{i}$ as the configuration of $\Psi_\L$ such that $\o_{i}|_{\L\setminus \obj(\ev)}=  {\o_{i-1}}|_{\L\setminus  \obj(\ev)}$ and
%${\o_{i}}|_{ \obj(\ev)}= \o'$
%}
%\vskip.2cm

\vskip.4cm

\fbox{
\begin{minipage}{0.9\linewidth}
{\textsc{Resampling}}.
\vskip.2cm
\begin{enumerate}
  \item Take a random evaluation $\o_0\in \Psi_\L$.
  \item  While there is a bad event belonging to $\fa$ occurring, select an event $\ev$  and resample all variables $\{\psi_x\}_{x \in \ \obj(\ev)}$.
  \item End while.
  \item Output current evaluation.
\end{enumerate}
\end{minipage}}

\vskip.4cm

\\Moser and Tardos proved
that if condition (\ref{cond}) of  Theorem \ref{lov} holds, then \textsc{Resampling}   terminates rapidly finding
a configuration $\o\in \Psi_\L$ such that none of the bad events of the family $\fa$ occurs. Later, inspired by the paper by Bissacot et al. \cite{BFPS},  Pegden \cite{Pe}
improved the Moser-Tardos result  replacing
condition (\ref{cond}) with condition (\ref{cond1n}).

\begin{theo}[Pegden]\label{MosTar} Given  a finite set $\L$  and its associated family
of mutually independent random variables $\psi_\L$,
let  $\FF$ be a family of tempered bad events  with natural dependency graph
$\GG$.
Let  $\bm{\mu}=\{\mu_\ev\}_{\ev\in \fa}$ be  non-negative numbers. If, for each $\ev\in \fa$,
\begin{equation}
\label{cond1np}
\prob(\ev)\;\leq\;  {\mu_\ev \over  \Xi_\ev(\bm \m,\GG)}
\end{equation}
with $ \Xi_\ev(\bm \m,\GG)$ defined in (\ref{Xi2}), then
$$
{\prob}\Big(\bigcap_{\ev\in \fa}\bar \ev\Big)\,> \,0
$$
and
algorithm \textsc{Resamplig} finds $\o\in\bigcap_{\ev\in \fa}\bar \ev$ in an  expected total number of  steps less than or equal to
$\sum_{\ev\in \FF}\m_\ev$.
\end{theo}
\\It is also worth to mention  that Kolipaka and Szegedy  showed  in \cite{KS2}  that  algorithm \textsc{Resampling} is successful in polynomial time also  if Shearer conditions hold (see also \cite{AP} for a similar result).
\vskip.2cm
\\{\bf Remark}. In the Moser-Tardos setting
above described the function $\Xi_\ev(\bm \m,\GG)$, defined in (\ref{Xi2}), admits
a somehow natural  upper bound so that condition (\ref{cond1np}) can be greatly simplified (paying some price as we will see in  a moment).
Just observe that, for  any $y\in {\L}$, the set   $\FF(y)=\{\ev\in \FF:~y\in \obj(\ev)\}$
is a clique of the natural dependency graph $\GG$ of the family $\fa$: any pair $\{\ev,\ev'\}\subset \FF(y)$ is such that ${\rm supp}( \ev)\cap {\rm supp}( \ev')\supset \{y\}\neq\0$, i.e., any pair $\{\ev,\ev'\}\subset \FF(y)$  is an  edge of $\GG$.
Thus, the neighborhood $\G^*_\GG(\ev)$ of any event  $\ev\in \FF$ is  the union of cliques $\{\FF(y)\}_{y\in \obj(\ev)}$, and then we can write
$\G^*_\GG(\ev)= \cup_{y\in \obj(\ev)}\FF(y)$.
Therefore
\be\label{clique}
\Xi_\ev(\bm \m,\GG)\le  \Xi^{\rm clique}_\ev(\bm \m,\GG)\equiv\prod_{y\in \obj(\ev)}\Big[1+ \sum_{\ev'\in \fa(y)}\m_{\ev'}\Big].
\ee
%A sligtly better clique estimate is
%\be\label{cliqueb}
%\Xi_\ev(\bm \m,\GG)\le \widetilde\Xi^{\rm clique}_\ev(\bm \m,\GG)\equiv\m_\ev +\prod_{y\in \obj(\ev)}\Big[1+ \sum_{\ev'\in \fa(y)\atop \ev'\neq \ev}\m_{\ev'}\Big]
%\ee

\\Using the expression $\Xi^{\rm clique}_\ev(\bm \m,\GG)$, which is much simpler to evaluate, in place of $\Xi_\ev(\bm \m,\GG)$ in condition (\ref{cond1n}),   it was possible to improve estimate for latin transversal
in \cite{BFPS} and  improve bounds
for  several chromatic indices in  \cite{NPS} and \cite{BKP}. It is important however to stress that
the estimate  (\ref{clique}) is efficient only if the cliques $\{\fa(y)\}_{y\in \obj(\ev)}$ are not too overlapped. In general this union is not disjoint, but when cliques  $\{\fa(y)\}_{y\in \obj(\ev)}$  are too overlapped, the replacement $\Xi_\ev(\bm \m,\GG)$ with $\Xi^{\rm clique}_\ev(\bm \m,\GG)$   tends to be a too crude estimate. This situation occurs for example
in the case of perfect and separating hash families (see \cite{PS}).

\vskip.2cm

\subsection{\bf The subset gas}
In the context of the connection between the LLL and the APS it is natural to ask  where the Moser-Tardos setting defined above does fit.
We explain below that the APS counterpart of the Moser-Tardos setting is the so-called  subset gas.

\\The subset gas, originally proposed by Gruber and Kunz \cite{GK},   is a particular realization of the abstract polymer gas which
appears in many physical situations and it is defined as follows.

\\Given  a countable set $\VU$, the space of polymers $\PP$ is defined as the collection of all finite subsets of $\VU$, namely
$\PP=\{\g\subset \VU: |\g|<+\infty\}$. The Boltzmann factor is then defined as $W(\g,\g')=0$ if $\g\cap\g' \neq\0$ and $W(\g,\g')=1$ if $\g\cap\g'=\0$. Thus, in such a  realization of the APS,  polymers have
a cardinality, so that one can speak about large  polymers and small polymers. Of course, as before, each
polymer $\g$   has an associated  activity $w_\g$.
In most of the physical realizations, $\VU$ is the vertex set of
a (possibly  infinite)
graph, typically the cubic lattice $\Z^d$ with the edge set being the set of nearest neighbor in $\Z^d$.
Actually, the subset gases appearing in the framework of statistical mechanics and specifically spin systems on $\Z^d$
have in general   a further characteristic. Typically,  in each site $x\in  \Z^d$  is defined a random variable $s_x$ (the {\it spin} at $x$)
taking values in some space $S_x$ and frequently this space  $S_x$ is the same for all $x\in \VU$. Then  the space of the polymers  $\PP$ is formed by the  pairs
$\g = ({\rm supp}(\g) , s_\g)$, where ${\rm supp}(\g)$ is, as before,
a finite subset of $\VU$ and
$s_\g$ is the  {\it spin configuration} of $\g$, i.e.  a function from ${\rm supp}(\g)$ to $\prod_{x\in {\rm supp}(\g)}S_x$. Classical  examples are the “thick” contours of
Pirogov-Sinai theory (see e.g. \cite{Si}, Chap. II]).
%Typically the activity of a colored polymer $\g$ is proportional to the factor $k^{-|{\rm supp}(\g)|}$.

\\The reader can see at this point the evident parallel with the variable setting of the LLL. Namely,
polymers $\g = ({\rm supp}(\g) , s_\g)$ in the subset gas correspond to (elementary) events  $\ev$ in the Moser-Tardos variable setting.

\\Generally, as far as the subset gas is concerned, the bound (\ref{clique}) is always used.
So, for the subset gas,  the condition (\ref{fp}) can be  written as
\be\label{fpgas}
 |w_\g| \leq \frac{\mu_\g}{\prod_{x\in \obj(\g)}\Big[1+ \sum\limits_{\g'\in \PP\atop x\in \obj(\g')}
 \m_{\g'}\Big]},~~~~~~~~~~\forall \g\in \PP.
 \ee
 Note that the above criterion is constituted by many  inequalities, i.e. as many inequalities as the number
 of total polymers, so if $\PP$ is infinite this number can be infinite. However, in the specific case of the subset gas
 this set of inequalities can be replaced (and usually is!) by a unique global inequality to which the set of activities
 must obey. Indeed, since by (\ref{fpgas}) we necessarily have that $\m_\g> |w_\g|$, a typical choice  is to set
 $$
 \m_\g= |w_\g|e^{a|\obj(\g)|}
 $$ with $a>0$. Such a choice permits to resume
the set of conditions (\ref{fpgas}) in terms of  a simple ``global" conditions on the set
of activities $\bm w$. Namely, the thesis of Theorem \ref{bis} holds, if
\be\label{grub}
\sup_{x\in \VU}\sum_{\g\in \PP\atop x\in \obj (\g)}|w_{\g}|e^{a|\obj(\g)|}\le e^a-1 ~~~~\mbox{for some $a>0$.}
\ee
%This global condition can be further simplified by setting
%$$
%E_\PP=\{n\in \mathbb{N}: \exists \g\in \PP~ s.t. ~  |\obj \g|=n\},  ~~~\PP^x_n=\{\g\in \PP: x\in \obj(\g)~{\rm and}~|\obj(\g)|=n\}
%$$
%$$
%w_n= \sup_{x\in \mathbb{V}}\sup_{\g\in \PP^x_n} |w_\g|~~~~~~~d_n=\sup_{x\in \PP}|\PP^x_n|
%$$
%According with these notations (\ref{grub}) is satisfied if
%\be
%\sum_{n\in E_\PP}w_n d_n e^{an}\le e^a-1 ~~~~\mbox{for some $a>0$}
%\ee
\\The above discussion on the connection between the LLL and the APS leads to conclude that
 the set of conditions of the  CELL criterion can be reexpressed in terms
of a global unique condition of the probabilities $\prob(\ev)$ of the events $\ev$ as far as we are in the Moser-Tardos variable setting.

\begin{lem}\label{global}
Given  a finite set $\L$  and  a family
of mutually independent random variables $\psi_\L$,
 let $\FF$ be a family of tempered events  and for $\ev\in \FF$ let
$\prob(\ev)$ be its probability in the product space generated by the variables $\psi_\L$.
If it is possible to find $a>0$ such that
\be\label{grub3}
{\sup_{x\in \L}\sum_{\ev\in \FF\atop x\in \obj(\ev)}\prob(\ev)e^{a|\obj(\ev)|}}\le e^{a}-1,
\ee
then
$$
{\prob}\Big(\bigcap_{\ev\in \fa}\bar \ev\Big)\,> \,0
$$
and
algorithm \textsc{Resampling} finds a configuration $\o\in\bigcap_{\ev\in \fa}\bar \ev$ in an  expected total number of  steps less than or equal to
$\sum_{\ev\in \FF}\prob(\ev)e^{a|\obj(\ev)|}$.
 \end{lem}
% Again, by introducing the following notations
%\be\label{nota}
%p_n=\max_{\ev\in F\atop|\obj(\ev)|=n}\prob(\ev), ~~~d_n=\max_{x\in \L}|\{\ev\in F:~|\obj(\ev)|=n~{\rm and}~x\in \obj(\ev)\}
%\ee
%\be
%E_F=\{n\in \mathbb{N}: \exists \ev\in F~{\rm s.t.}~|\obj(\ev)|=n\}
%\ee
%condition (\ref{grub3}) is satisfied if
%\be\label{grub4}
%\sum_{n\in E_F}p_nd_n e^{an}\le e^a-1
%\ee

\\The latter  global ``subset gas condition" (\ref{grub3}) is able to reproduce  all improvements obtained via CELL
(\cite{NPS, BKP, BFPS}) with the exception of \cite{PS}, where the general condition  (\ref{cond1n}) has been used. It is worth to mention that
similar (but less effective) global conditions deduced from the original LLL  in the variable setting  have been already formulated in the literature and used in specific examples (see  e.g. Lemma 3 in  \cite{KS1} and reference therein).

\vskip.2cm
\subsection{The entropy compression method}

\\The Moser-Tardos algorithmic version of the LLL, since its appareance, has been the subject of a very intense study by several researchers
in the areas of computer science, combinatorics and probability. In this regard two main directions can be pointed out. The first one concerns
the (succesful) efforts  made to extend the validity of the algorithmic version of the LLL beyond the variable setting
(see e.g. \cite{HV}, \cite{AI0}, \cite{HSi}, \cite{AIK},  \cite{Ko}, \cite{Har} and references therein)
with the objective to include
important applications of the non constructive LLL (such as latin transversal) which does not fit in the independent variable setting.

\\The second direction was motivated by the  fact that the algorithm \textsc{Resampling} proposed by Moser and Tardos was extremely  simple, so it makes sense to try to modify/refine  it in order to improve the  final criterion beyond   (\ref{cond1n}).  These ideas have been
originally developed in \cite{DJKW, GKM, EP}  where   backtracking algorithms
have been
implemented  for  specific graph coloring problems to obtain
bounds which are  better than those  obtainable  by   LLL or   CELL.

\\In particular,  Esperet and
Parreau  devised in \cite{EP} an algorithm able to obtain
a new upper bound for the chromatic index  of the acyclic edge coloring  of a graph with maximum degree $\D$ which sensibly
improves on the bound obtained by Ndreca et al.  \cite{NPS}  just an year before via the CELL.
The algorithm proposed by Esperet and Parreau  presents evident differences  from the
algorithm \textsc{Resampling}: instead of sampling all variables at once  and then resampling some of them  until all flaws are avoided, the Esperet-Parreau algorithm  starts with the  empty configuration in  which no variable has an assigned value. Then, step by step
a  (random) value  is attributed  to a currently unassigned variable; when  this leads to the appearance of one or more flaws, the algorithm backtracks to a partial non-violating configuration by retracting some set of variables.

\\Esperet and
Parreau
suggested, through further examples and applications, that their  algorithm   could be adapted  to treat most of the
applications in graph coloring problems covered by the LLL.
Indeed, this was confirmed in several successive  papers   \cite{GMP,Pr,OP,PSS,CR, CLY, SV, FrS, BBCFGM}, where the
Esperet-Parreau scheme has been applied to various
graph coloring problems and beyond, generally improving previous results obtained via the LLL/CELL
(sometimes the improvement is more sensible, sometimes less).
However, in all  papers mentioned above the Esperet-Parreau algorithmic scheme, usually called {\it entropy compression  method} (the name is probably due to Tao \cite{Ta}), has been  commonly
 utilized as  a set of {\it ad hoc} instructions  to be implemented on a case-by-case basis. A systematization of the entropy compression method
 providing a general criterion similar to those given by the LLL and CELL has been an open question since the beginning (\cite{EP}, \cite{GMP}) and it is still demanded even  in the very recent paper by  Achlioptas and Iliopoulos \cite{AI}. In this regard we mention
a non algorithmic general criterion proposed by Bernshteyn \cite{Be} which is able to reproduce several results obtained by the entropy compression method. The systematization of the Esperet-Parreau backtrack algorithm
can actually be found   in  \cite{APS} in which the setting where entropy compression can be used is clearly outlined,
a general
{\it entropy compression criterion} is proposed and  a  connection  between bounds obtained via this method and LLL conditions  is elucidated.

\\According to \cite{APS},  the entropy compression method can be implemented  in any application that can be analyzed through   a variable setting {\it a la} Moser-Tardos  with the further restriction that, for all $x\in \L$,
variables $\psi_x$  take values in a common space $\Psi_x=[k]$, where $k\in \mathbb{N}$
(or, more generally, in possibly distinct spaces $\Psi_x$ but all with the same cardinality $k$)   according to the uniform distribution. We will refer to this particular realization of the variable setting as
the {\it uniform variable setting}. Let us review rapidly for later comparison the
entropy compression criterion obtained in \cite{APS}.

\subsubsection{The entropy compression  setting}
As said above, the entropy compression method can be applied in the so-called {\it uniform variable setting} where all random variables $\{\psi_x\}_{x\in \L}$ take values in the common space $\Psi=[k]\equiv \{1,2,\dots, k\}$ according to the uniform distribution.   We sometimes  refer  to $[k]$ as the set of ``colors" and we set $[k]_0=[k]\cup\{0\}$. So the entropy compression setting is determined by the pair $(\L, k)$, where $\L$ is a finite set and $k$ is a positive integer.
A configuration $w$  of $X\subset \L$ is a function
$w:X\to [k]$, and a partial configuration $w$ of $X\subset \L$ is a function $w:X\to [k]_0$ and
when  $w(x)=0$ we say that the variable $\psi_x$ is unassigned (or uncolored).
For any non empty $X\subseteq \L$, let $[k]^X$ and $[k]_0^X$ denote
the sets of configuration and partial configuration
in $X$ respectively. Given $Y\subset X$ and $w\in [k]^X$ we denote by $w|_{Y}$ the restriction
of the configuration $w$ to $Y$. Of course $w|_{Y} \in [k]^Y$.
\vskip.2cm
\\{\bf Remark}. More generally, one can also suppose that variables   $\{\psi_x\}_{x\in \L}$ take values in possibly different spaces $\{\Psi_x\}_{x\in \L}$ but all having a common cardinality $k$, i.e. such that $|\Psi_x|=k$ for all $x\in \L$ and in each space $\Psi_x$ the random variable $\psi_x$ takes values according to the uniform distribution. The Example 2 in Section \ref{appl} below falls in this more general setting.
\vskip.2cm
\\Given the pair $(\L,k)$ and $A \varsubsetneq \L$, a {\it  flaw of $A$}  is a subset $\ev\subset [k]^A$, i.e, it is a tempered  bad event in the Moser-Tardos scheme. As before,  we write $\obj(\ev)=A$ and we call $|\obj(\ev)|$ {\it the size of $\ev$}
while $|\ev|$ is the number of configurations forming $\ev$. An event is {\it elementary} if $|\ev|=1$. As usual, $\bar \ev$ will denote the complement of $\ev$ in $[k]^A$,
i.e. $\bar\ev=[k]^A\setminus\ev$.

\\Given  a family $\fa$ of tempered events,
a {\it a good configuration w.r.t. to $\fa$} is a configuration $\o\in \O\equiv[k]^\L$  avoiding all events in $\fa$, that is to say,
$\o$ is such that for  all $\ev\in \fa$ such that $\obj(\ev)=A$, we have $\o|_A\notin \ev$. Given $\ev\in \FF$ and $X\varsubsetneq \obj(\ev)$
we denote  $\ev|_X=\{w_X\in [k]^X: w\in \ev\}$.

\vv
\\{\bf Remark}.
We stress once again that virtually all
applications of the entropy compression method available  in the literature fall in this
uniform  variable setting, with the sole exception, as far as we know, of the
the acyclic edge coloring of  bounded degree graphs (see \cite{APS}).

\vskip.2cm

\begin{defi}\label{seeda}
Given an event $\ev$, a non empty subset $X\subset \obj(\ev)$ is called a seed of $\ev$ if $\ev|_X=[k]^X$ and $|\ev|=k^{|X|}$.
An event  $\ev$
is said tidy  if either it is elementary, or it is such that for all $y\in \obj(\ev)$,  there exists  a non-empty set   $X\subset \{\obj(\ev)\setminus \{y\}\}$ which is a seed of $\ev$.
\end{defi}
Clearly, by definition,  all seeds of a tidy  event $\ev$  must have  all the same cardinality  which we denote by $\k(\ev)$. If $\ev$ is elementary, we set $\k(\ev)=0$.
We further set
\be\label{power}
\|\ev\|= |\obj(\ev)|-\k(\ev)
\ee
and refer to $\|\ev\|$ as {\it the  power of the event $\ev$}.
Note that  $\|\ev\|=|\obj(\ev)|$ if and  only if
 $\ev$ is elementary. Moreover, the properties listed here below follow immediately from Definition \ref{seeda}.
\begin{enumerate}

\item  If $\ev$ is tidy and  $X$ is a seed of $\ev$, then any configuration $w\in \ev$ is uniquely determined by its restriction to $X$
and no $Y\varsubsetneq X$ has this property.

\item If $\ev$ is tidy, then
$$
\prob(\ev)={1\over k^{|\obj(\ev)|-\k(\ev)}}.
$$

\item Let $\FF$ be a family of tempered and tidy events. Suppose that  an event $\ev$ is occurring in a given configuration and $X\varsubsetneq\obj(\ev)$ is a seed of $\ev$.  If we resample all the variables in $\obj(\ev) \setminus X$ leaving all other variables unchanged,  then in the so obtained new configuration
    the probability for an event $\ev' \in \FF$ to occur is  $\prob(\ev')$.
\end{enumerate}
\\{\bf Remark 1}. In \cite{APS} the definition of seed is slightly more general. Namely it {coincides with item 1} of the above list:
a seed  $X\subset \{\obj(\ev)\setminus \{y\}\}$ is
such that any coloring $w\in \ev$ is uniquely determined by its restriction to $X$
 and no $Y\varsubsetneq X$ has this property.  The difference is subtle. As an example, suppose that the edges of a graph $G=(V,E)$ are
colored  at random  using $k$ colors, uniformly  and independently and consider a cycle $C$  of $G$ with an even number of edges and then let the event $\ev_C$ be ``$C$ is properly bichromatic'', i.e. its edges are colored with two colors and no pair of adjacent edges are monochromatic.
 Then any cherry (i.e. two incident edges)  of $C$ is a  seed of $\ev_C$ according to the definition given in \cite{APS}, while $\ev_C$ is not tidy
 according to the definition \ref{seeda}, since the restriction of $\ev_C$ to any cherry  $c\in C$ is not in $[k]^c$, once the
 monochromatic configurations of $c$ are not allowed if $C$ is properly bichromatic.
 On the other hand given a path $p$ of $G$ constituted by an even number of edges, let $\ev_p$ be the event ``the second half of $p$ is colored in the same way as the first half'', then $\ev_p$ is tidy according to the Definition \ref{seeda}. For example, we can set as a seed the set of edges constituting the first half of the path or the set of edges constituting the second  half of the path.

 \vskip.2cm

\\{\bf Remark 2}.
 By definition, a tidy flaw has  the empty set as its unique seed if and only if is elementary. Note that
if $\ev$ is not tidy, then it can be seen as  the disjoint union of tidy (in the worst case  elementary) flaws.
Therefore there is no loss of generality in considering only
families  in which all flaws are tidy.

\vv
\\We introduce the following notations.
$$
\fa_{s}(y)=\{ \ev\in \fa:~y\in \obj(\ev)~{\rm and}~\|\ev\|=s\}
$$
and
\be\label{dj}
d_s= \max_{y \in \L} \big|\fa_{s}(y)\big|.
\ee
Namely, $d_s$ is an upper bound for the number of events with power $s$ whose support  contains a common element of $\L$.
\\We finally define
\be\label{EA}
E'_\FF=\{s\in \mathbb{N}: \exists \ev\in \fa~  {\rm such ~that }~  \|\ev\|=s\}.
\ee
%and set,
%for $\a\in (0,+\infty)$,
%\be\label{phiea}
%\varphi_F(\a)=  1+ \sum_{j\in E'_F}d_j \a^{j}
%\ee

\vv\vv

%We finally remark that, once introduced the definition of $\k_l$, we have, for free,
%a upper bound for the probability of each event $A$, i.e. if $A\in {{\mathcal O}}^{\k_l}_l$, then
%$P(A)\leq {k^{-(l-\k_l)}}$.

\noindent
%Let ${\mathcal O}$ be a set of events depending on variables $\cal V$.

\subsubsection{The entropy compression algorithm and the entropy compression Lemma}

\\We assume  that a total order has been chosen in the sets  $\L$  and $\fa$.
We  choose, for each $y\in \L$ and $\ev\in \fa(y)$,  a
unique subset  $G(\ev, y)\subset \obj(\ev)\setminus\{y\}$  such that  $G(\ev,y)$ is a
seed of $\ev$. We also denote
shortly $G^{c}(\ev,y)=\obj(\ev)\setminus G(\ev,y)$. Note that $y\in G^c(\ev,y)$.
Given a partial coloring $w\in [k]_0^{\L}$, given $X\subset \L$
  and given  $s\in [k]\cup\{0\}$ we denote by $w]^s_{_X}$
  the partial coloring which coincides with $w$ in the set $\L\setminus X$
  and it takes the value
$s$ at every   $x\in X$. If $X=\{x\}$ we set shortly $\omega]^s_{X}\equiv \omega]^s_{{x}} $.

\vskip.2cm
\noindent
Let $t$ be an arbitrary natural number (which can be taken as large as we please)
and let   $V_t$ be an element
of  $[k]^t$, i.e. $V_t$ is a vector with $t$ entries such that each entry takes values in the set $[k]$.

\\The algorithm \textsc{entropy compression} has input $V_t$,
performs (at most) $t$ steps and, at  each step $i\in [t]$,  produces a partial coloring $w_i$  as described below.

\vskip.2cm

\fbox{
\begin{minipage}{0.9\linewidth}
\textsc{entropy compression} (with input $V_t$)
\begin{itemize}
  \item [-] Step $0$. Set $w_0=0$, i.e. in the beginning no element $\psi_x$, with $x \in \L$, is colored.
  \item [-] Step $i$ (for $i\ge 1$).
  \item [{$i^\circ$})]
  If $w^{-1}_{i-1}(0)\neq\emptyset$,
  let $y$ be  the smallest  element of $\L$ (in the total order chosen) such that $\psi_y$ is uncolored in the partial coloring $w_{i-1}$.
  Take the $i^{th}$ entry of the vector $V_t$ and  let $s\in [k]$ be this entry. Color
  $\psi_y$ with the  color $s$ and consider the partial coloring $w_{i-1}]^s_y$ obtained from $w_{i-1}$ by coloring $\psi_y$ with the color $s$.
  \begin{enumerate}
  \item[$i^\circ_1)$]
  If no flaw  occurs in $w_{i-1}]^s_y$, set
  $\omega_i=\omega_{i-1}]^s_y$ and go to the step $i+1$.
  \item[$i^\circ_2)$]
  Conversely, if some  flaw occur in $w_{i-1}]^s_y$, select the smallest, say $\ev$, which by construction belongs to the set $\fa(y)$. Set $w_i=w_{i-1}]^0_{G^c(\ev,y)}$ and go to the  step $i+1$. In words, $w_i$ is obtained from $w_{i-1}$ by  discoloring all $\psi_x$ such that $ x \in \obj(\ev) \setminus G(\ev, y)$.
\end{enumerate}
\item[{$i^\bullet$)}]
If $w^{-1}_{i-1}(0)=\emptyset$, stop the algorithm discarding all entries $v_i,v_{i+1},\dots,v_t$ of $V_t$.
\end{itemize}
\end{minipage}
}

%\\\textsc{entropy compression} (with input $V_t$)
%\vskip.1cm
%\\- Step $0$. Set $w_0=0$, i.e. in the beginning no element $\psi_x$, with $x \in \L$, is colored.
%\\- Step $i$ (for $i\ge 1$).
%\begin{itemize}
%\item [{$i^\circ$}]
%If $w^{-1}_{i-1}(0)\neq\emptyset$,
%let $y$ be  the smallest  element of $\L$ (in the total order chosen) left uncolored by the partial coloring $w_{i-1}$.
%Take the $i^{th}$ entry of the vector $V_t$ and  let $s\in [k]$ be this entry. Color
%$y$ with the  color $s$ and consider the partial coloring $w_{i-1}]^s_y$ obtained from $w_{i-1}$ by coloring $y$ with the color $s$.
%\begin{enumerate}
%\item[$i^\circ_1)$]
%If no flaw  occurs in $w_{i-1}]^s_y$, set
% $\omega_i=\omega_{i-1}]^s_y$ and go to the step $i+1$.
%\item[$i^\circ_2)$]
%Conversely, if
%some  flaw occur in $w_{i-1}]^s_y$,
%select the smallest, say $\ev$, which by construction
%belongs to the set $\fa(y)$.
%Set $w_i=w_{i-1}]^0_{G^c(\ev,y)}$ and go to the  step $i+1$.
%In words, $w_i$ is obtained from $w_{i-1}$ by  discoloring all elements in $\obj(\ev)$
%except those belonging to the set $G(\ev, y)$.
%\end{enumerate}
%\item[{$i^\bullet$}]
%If $w^{-1}_{i-1}(0)=\emptyset$, stop the algorithm discarding all entries $v_i,v_{i+1},\dots,v_t$ of $V_t$.
%\end{itemize}

\vskip.4cm

\\Note that the partial coloring $w_i$  returned by the algorithm at the end of each step $i$  necessarily avoids all flaws in $\fa$.
\textsc{entropy compression} performs  at most  $t$ steps but it can stop  earlier, i.e. after having performed $m<t$ steps
and $\omega_m^{-1}(0)=\emptyset$.
In this case only the first $m$ entries of the vector $V_t$ are used. \textsc{entropy compression} is successful if it stops after $m<t$ steps,
or it lasts $t$ steps and after the last step $t$ we have $\omega_t^{-1}(0)=\emptyset$. Conversely,
\textsc{entropy compression} fails if it performs all $t$ steps and $\omega_t^{-1}(0)\neq \0$. Clearly when \textsc{entropy compression} is successful $w_t$
is a good configuration.
Observe that  \textsc{entropy compression} can be either deterministic, if $V_t$  is a given prefixed vector, or random,
if the entries of $V_t$ are uniformly sampled from the set $[k]$  sequentially and  independently. In \cite{APS} the following theorem is proved.

 \begin{theo}[{\bf Entropy compression  lemma}]\label{ecl}
  Assume that a  pair $(\L,k)$ is given together with
a family  $\fa$  of tempered and tidy flaws. If there is $\a>0$ such that
\begin{equation}\label{condentr}
{1+ \sum_{s\in E'_\fa} d_s\a^s\over \a}< k,
\end{equation}

\\then
$$
\bigcap_{\ev\in \fa}\bar \ev\neq\0.
$$
Moreover, \textsc{entropy compression} finds a configuration $\o\in \bigcap_{\ev\in \fa}\bar \ev$  in an expected number of steps linear in $|\L|$.
\end{theo}

\\We will refer to the inequality (\ref{condentr}) as the {\it entropy compression criterion}.  We stress that this theorem is able to reproduce
all results obtained in these last years via the entropy compression method.

\\It is now simple to compare the above entropy compression criterion with the global CELL criterion (\ref{grub3}). Since we are in the entropy compression setting determined by the pair $(\L,k)$, the probability space is generated
by $|\L|$ i.i.d. uniformly distributed random variables taking values in $[k]$.
As recalled above, in the restricted variable setting covered by the entropy compression method, we have that $\prob(\ev)={1\over k^{\|\ev\|}}\equiv p_{\|\ev\|}$. Setting
\be\label{qu}
q=\max_{\ev\in \fa}\left\{ {|\obj(\ev)|\over \|\ev\|}\right\}
\ee
and grouping events in terms of their powers we get
$$
\sup_{x\in \L}\sum_{\ev\in \FF\atop x\in \obj(\ev)}\prob(\ev)e^{a|\obj(\ev)|}\le {{\sum_{s\in E_\FF'} d_s p_s{e^{aqs}}}}.
$$
So condition (\ref{grub3}) is fulfilled  if
\be\label{nps2}
{{\sum_{s\in E'_\FF} d_s p_s{e^{aqs}}}}\le e^a-1,
\ee
or, setting $\a={e^{aq}\over k}$,
if there is $\a>0$ such that
\be\label{nps}
{(1+ \sum_{s\in E'_\FF} d_s\a^s)^q\over \a}\le k.
\ee
\\The reader can immediately compare (\ref{nps})  with  the entropy compression  condition (\ref{condentr}).
The presence of the exponent $q$ defined in (\ref{qu}) in inequality (\ref{nps}) is the only reason why entropy compression condition (\ref{condentr})  can  give better bounds than LLL. It must however be stressed that we are  excluding  here  the  case of the acyclic edge coloring of a graph with maximal degree $\D$. In this pretty singular case,  the entropy compression scheme is fruitfully combined with
the crucial observation that it is possible to properly color the edges of $G$ using just $2\D-1$ color in such a way to avoid
bichromatic cycles of length 4  (see Lemma 4 in \cite{APS}). This  fruitful strategy leaded Esperet and Parreau
to a very sensible   improvement of the upper bound of the acyclic edge chromatic index of $G$ with respect to the bound obtained via CELL.
Due to its specificity, the  case of the acyclic edge coloring  must be treated separately (see comments below and see also the remark in Section 4.2.5 of \cite{APS}).

\subsection{Motivations and plan of this paper}

\\Concluding this introduction, we  need to mention   two recent papers, \cite{KL} and \cite{GKPT}, proposing a variant of the Moser-Tardos resampling algorithm, which has  actually motivated  the present paper. In particular, in \cite{GKPT} Giotis et al.  are able to slightly improve just the specific case of the acyclic edge chromatic index of a graph with maximum degree $\D$.  The
intriguing fact is that Giotis et al. use in \cite{GKPT}  the Moser-Tardos resampling algorithm with the unique variant that the successive resampled bad events must be chosen, when possible, in the neighbor of the previous bad event.
Their result is somehow surprising considering that, as mentioned above,
the  CELL  criterion applied to acyclic edge coloring gives a much worse bound than entropy compression method.

\\In the present  paper we
manage to combine the  ideas of \cite{GKPT} (also foreshadowed  in \cite{KL})
with the observation explained above that the power of an event (possibly tidy) can be considered in place of the cardinality of its support
and we  show that
the  criterion  (\ref{condentr}) based on the backtracking algorithm \textsc{entropy compression} can be reobtained in the usual
 Moser-Tardos scheme by doing the slight modification of the  algorithm \textsc{Resampling}  illustrated  in \cite{KL} and \cite{GKPT}  (which give rise to forests instead of tree as a register of the steps of the algorithm) jointly with the prescription proposed by Esperet-Parreau to not resample certain variables of the bad events (the previously seen  ``seeds" of the events).

 \\In the very specific and  singular case of acyclic edge coloring,
 the Moser-Tardos modified algorithm presented in this paper is able to further slightly improve the  bound obtained in \cite{GKPT}, and this latter issue is the subject of a separate paper \cite{FLP}.

\\The rest of the paper is organized as follows. In Section \ref{vmt} we describe the variant of the Moser-Tardos algorithm and state our
main result, i.e. Theorem \ref{teo1}. Section \ref{proof} is devoted to the proof of Theorem \ref{teo1}.
Finally in Section \ref{appl} we present some examples.

\section{A variant of the Moser-Tardos Algorithm }\label{vmt}
\zeq

\\Let us consider the general Moser-Tardos framework. Given a finite set $\L$ with cardinality $m\equiv|\L|$, let $\psi_\L\equiv \{\psi_x\}_{x\in \L}$ be a set  of  $m$ mutually independent random variables such that each $\psi_x$ takes values
in $\Psi_x$ and let $\O_\L=\prod_{x\in \L}\Psi_x$
 be the product probability space generated by these variables. Let $\FF$ be a finite collection of tempered events in $\O$. We recall that for each event $\ev \in \FF$ there exists a subset $\obj(\ev) \varsubsetneq \L$   such that $\ev$ depends on variables $\{\psi_x\}_{x\in \obj(\ev)}$. As usual,  if $U\subset \L$, we set $\psi_U=\{\psi_x\}_{x\in U}$ and
 $\O_U= \prod_{x\in U}\Psi_x$. Moreover,
given a random configuration $\o\in \O_\L$,  $\prob(\ev)$ denotes  the probability
 of the event $\ev\in \FF$ to occur.  We also recall that if $\ev$ is a tempered event such that there is a unique configuration  $\omega \in \Om_{\obj(\ev)}$ which realizes $\ev$ we say that $\ev$ is an {\it elementary} event.

\begin{defi}[Seed]\label{seed}
Let $\ev\in \FF$. A proper subset $U \varsubsetneq \obj(\ev)$ is called a ``seed" of $\ev$ if  it is  such that,  given
a configuration $\o\in \O$ such that the event $\ev$ is occurring, if we resample all the variables in $\obj(\ev) \setminus U$ leaving unchanged the values of all the other variables, then in the new configuration
$\o'$ so obtained all the events $\ev' \in \FF$ have at most $\prob(\ev')$ to happen and any $U'\supset U$ has not this property. We denote by $S_{\ev}$ the set of all seeds of $\ev$.
An event $\ev\in \FF$ is tidy if it is such that all seeds of $\ev$ have the same non zero cardinality $\k(\ev)$ and for all $x\in \obj(\ev)$ there exists $U\in S_\ev$ such that $x\notin U$.
%any configuration of the values for the variables in $S(A)$ is still possible and any $U\supset S(A)$ has not this property.
\end{defi}

\\{\bf Remark}\label{od}.
 Of course,  in the uniform variable setting, where the variables $\{\psi_x\}_{x\in \L}$ beside being independents are also identically and uniformly distributed and taking values in the common  set $[k]=\{1, \dots, k\}$,  Definition \ref{seed} and  Definition \ref{seeda} are equivalent. It is also important to stress  that
Definition \ref{seed} says that once we reach a configuration $\o$ in which  the event $\ev$ occurs and $U$ is a non-empty seed of $\ev$, passing to a new configuration $\o'$  obtained from $\o$ by resampling only variables
$\{\psi\}_{\obj(\ev)\setminus U}$  does not help any event to happen. The transition $\o\to\o'$ reminds the definition of {\it resampling oracle} given in \cite{HV}.

\vskip.2cm

\\We define the {\it power} of  the event  $\ev$ as the number
\be\label{power2}
\|\ev\| = \cases{|\obj(\ev)|-\k(\ev) &if $\ev$ is tidy,\cr\cr
|\obj(\ev)| & otherwise.
}
\ee

\\Moreover,  for any tidy event $\ev\in \FF$ and any   $x \in \obj(\ev)$, we fix a rule to choose uniquely a seed  $S_{x}(\ev)$ of  $\ev$ such that $x\notin S_{x}(\ev)$. If $\ev$ is either elementary or not tidy we set $S_{x}(\ev)= \emptyset$. Note that in any case
\be\label{power3}
|\obj(\ev)|-|S_x(\ev)|=\|\ev\|.
\ee

 %\\Given an event $A$ and a variable $\psi$, in general the seed $S_{\psi}(A)$ is not unique, observe that in Example \ref{P} $S'$ and $S''$ are candidates for $S_{v_k}(A_p)$.  Then, given an event $A$ and a variable $\psi \in vbl(A)$, we fix a rule such that the seed  $S_\psi(A)$ is uniquely determinate.

%\\However, not all events $A$ and variables $\psi \in vbl(A)$ are such that $S_{\psi}(A)$ exists. In fact, consider a path composed by three vertices $p=\{v_1, v_2, v_3\}$ and let $A_p$ be the event ``the path $p$ is such that $v_1$ and $v_3$ are monochromatic''. Observe that $S_{v_1}(A_p)=\{v_2, v_3\}$ and $S_{v_3}(A_p)=\{v_1, v_2\}$ but it does not exist $S_{v_2}(A_p)$.

%\begin{rema}[{\red To think about}]
 % \begin{enumerate}
  %  \item if $S$ and $S'$ are seeds of an event, then they have the same cardinality? Or should we consider only events with this property? Or should we always take the smallest seed in terms of cardinality? I don't know...
   % \item In this paper we will consider that given an event $A$, for all variable $\psi \in vbl(A)$ we can find a seed $S_{\psi}(A)$ such that $\psi \notin S_\psi(A)$. If the event $A$ does not have this property, we say also that $S=\0$ is the unique seed of $A$.
    %  \end{enumerate}
%\end{rema}

\\We will classify the events $\ev \in \FF$ according to their power  $\|\ev\|$. Let  ${{E'_\FF}}\subset\mathbb{N}$  be defined  as
\be\label{E'}
E'_{\FF}=\{s \in \mathbb{N}: \exists \ev\in  \FF ~s.t. ~\|\ev\|=s\}.
\ee
For  $s\in E'_\FF$, we set
%\begin{equation}
%{\cal A}^l_k=\{A\in {\cal A}: |A|=k~{\rm and}~\|A\|=l\},
%\end{equation}
\begin{equation}\label{s}
  \FF_s=\{\ev\in \FF: \|\ev\|=s\},
\end{equation}
%and for $k\ge 1$, let $L_k$ be the set of all possible cardinalities of seed of events of size $k$, i.e.
%\begin{equation}
%L_k=\{l\in \mathbb{N}: \exists A\in {\cal A}^l_k \}.
%\end{equation}
Finally, for  $x\in \L$ and $s\in E'_{\FF}$ , we define $d_s(x)$ as %the number of events $A \in {\cal A}^l_k$ such that $\psi\in vbl(A)$,
%\begin{equation}
%d^l_k(\psi)=|\{A\in {\cal A}^l_k: \psi\in vbl(A)\}|,
%\end{equation}
\begin{equation}
  d_s(x)=|\{\ev\in \FF_s: x\in\obj(\ev)\}|,
\end{equation}
and set
%\begin{equation}
%d^l_k=\max_{\psi\in \Psi} {d^l_k(\psi)}.
%\end{equation}
\begin{equation}
  d_s=\max_{x\in \L} {d_s(x)}.
\end{equation}
\def\0{\emptyset}
\noindent
%Let $\cal A$ be a set of events depending on variables $\cal V$.

\\Hereafter we will assume that  a total order is fixed in the set $\L$ as well as on the set of events $\FF$. %Given a subset $\Psi' \subset \Psi$, when $\psi$ is the variable with the smallest label in $\Psi'$, we will say shortly that $\psi$ is the smallest variable in $\Psi'$, and we also define the smallest event in a subset $\A' \subset \A$ as the event with the smallest label in $\A$.
%For each $A\in {\cal A}$, we choose a unique subset  $S_A\subset vbl(A)$  such that  $S_A$ is a seed of $A$ (if $A$ has no seeds, then $S_A=\0$).
\\Following \cite{GKPT},  we now describe a procedure, called \textsc{Forest-Algorithm} which  samples (and eventually resamples) the variables  $\psi_\L$. %Each discrete time $t\in \mathbb{N}$ such that \textsc{Forest-Algorithm} samples (or resamples) a variable is called an {\it instant}.
Given an evaluation $\o$ of all variables $\psi_\L$, we say shortly that the atom $x\in \L$ is {\it bad} if  some $\ev\in \FF$ occurs in the evaluation $\o$ and $x\in \obj(\ev)$. Otherwise we say that $x$ is {\it good}.

\vskip.5cm

\fbox{
\begin{minipage}{0.9\linewidth}
\vskip.2cm
{\textsc{ Forest-Algorithm}}.
\vskip.2cm
\begin{enumerate}
  \item Sample all variables  $\psi_\L$.
  \item  While there is a bad atom, select the pair $(x, \ev)$ where $x$ is  the smallest bad atom and where $\ev$ is the smallest event occurring such that $x \in \obj(\ev)$, and do
  \item  \textsc{Resample}$(x,\ev)$.
  \item End while.
  \item Output current evaluation.
\end{enumerate}
\end{minipage}
}
\vskip.6cm
\fbox{
\begin{minipage}{0.9\linewidth}
\vskip.2cm
{\textsc{ Resample}$(x,\ev)$}
\vskip.2cm
\begin{enumerate}
  \item Resample all variables $\psi_y$ such that $y \in \obj(\ev)\setminus S_x(\ev)$.
  \item While there is a bad atom in $\obj(\ev)\setminus S_x(\ev)$, let $x'$ be the smallest of these atoms  and  let $\ev'$ be the smallest event occurring such that $x' \in \obj(\ev')$ and do
  \item \textsc{Resample}$(x',\ev')$.
  \item End while.
\end{enumerate}
\end{minipage}
}
\vskip.5cm

%\begin{rema}
 % The resample of the variable $\psi$ is uniformly and independently among all the possibilities. Now we do not avoid ``colors''.
%\end{rema}

\\A {\it step} of \textsc{Forest-Algorithm}  is the procedure described in Line 2 of \textsc{Resample}$(x,\ev)$.
Observe that, since  $x\notin S_x(\ev)$ for any $\ev$ such that $x \in \obj(\ev)$, in \textsc{Resample}$(x,\ev)$ the variable $\psi_x$ is always resampled.
A {\it phase } of \textsc{Forest-Algorithm} is the collection of steps made by \textsc{Forest-Algorithm} during a  call of \textsc{Resample}$(x,\ev)$ in Line 3 of \textsc{Forest-Algorithm}. Note that during a phase many steps occur, the first step of the $i$-th phase will be called the {\it root} of the phase $i$. The {\it record} of the algorithm is the list $$\mathcal{L}=((x_1,\ev_1), (x_2,\ev_2),\dots )$$ constituted by the steps done by the algorithm during its execution. We will denote by {\it atom label} (resp. {\it event label}) any atom (resp. event) listed in the record $\mathcal{L}$. According to the prescriptions described above, $\mathcal{L}$ is a random variable determined by the random samplings performed by the algorithm in each step. If  $\mathcal{L}$ is finite, i.e. if $|\mathcal{L}|=n$ for some $n\in \mathbb{N}$,  then the algorithm terminates having performed $n$ steps and produces an evaluation $\o \in \bigcap_{\ev \in \FF}\bar\ev$. Let us define

\begin{equation}\label{p1}
  P_n= \prob(|\mathcal{L}|=n).
\end{equation}

\\In other words $P_n$ is the probability that \textsc{Forest-Algorithm} runs $n$ steps.

\\We are now in the position to state the main result of this paper. To do this we introduce the following notations. For $s\in E'_\FF$ and $\xi>0$, define
 \begin{equation}\label{pij}
 p_s=\max_{\ev\in \FF_s} \prob(\ev),
 \end{equation}

\be\label{fea}
 \phi_{\FF}(\xi)= \sum_{s\in E'_\FF}p_s d_s (\xi+1)^{s}.
\ee

 \begin{theo}\label{teo1}
 Given  a finite set $\L$  and  a family of mutually independent random variables $\psi_\L$,
 let $\FF$ be a family of tempered events depending on $\psi_\L$.
 Suppose that
\be\label{teo}
 \min_{\x>0} {\phi_{\FF}(\x)\over \x} <1,
\ee
then there is an evaluation of the variables $\psi_\L$ such that none of the events in the family $\FF$ occur. Moreover,  \textsc{Forest-Algorithm} finds a configuration $w\in \bigcap_{\ev\in \fa}\bar \ev$  in an expected number of steps polynomial  in $m=|\L|$.
\end{theo}

\\{\bf Remark}. Note that (\ref{teo}) is completely equivalent to the condition (\ref{condentr}) of entropy compression lemma (Theorem \ref{ecl}). Indeed, in the uniform variable setting defined by the pair $(\L,k)$, we have that $p_s={1\over k^s}$ and therefore, posing $\a={(\x+1)/ k}$, condition (\ref{teo})  is
rewritten
in the form (\ref{condentr}).

\section{Proof of Theorem \ref{teo1}}\label{proof}

Let us start by proving some important properties of \textsc{ Forest-Algorithm}.

\begin{lem}\label{Agood2}
Consider any call of \textsc{Resample}$(x,\ev)$ and  let $Y$ be the set of all good atoms at the beginning of this call. If this call finishes, then all the atoms  in $Y\cup \{\obj(\ev) \setminus S_x(\ev)$\} are good.
\end{lem}
{\it Proof:} According to the algorithm if \textsc{Resample}$(x,\ev)$ finishes then $\obj(\ev)\setminus S_x(\ev)$ are good atoms, so we just need to prove that the atoms  in $Y$ continue to be  good in the end of \textsc{Resample}$(x,\ev)$. Let $y\in  Y$, and assume that \textsc{Resample}$(x,\ev)$ finishes and performs $n$ steps. Suppose by contradiction that after these $n$ steps performed
by \textsc{Resample}$(x,\ev)$  $y$ is bad.
Then there exists a last step $t\le n$ of \textsc{Resample}$(x,\ev)$  such that $y$ was good at step $t-1$, became  bad  at step $t$ and stayed  bad during the remaining $n-t$ steps of \textsc{Resample}$(x,\ev)$.
This means that there is an event $\ev'$ and an atom $z\in \obj(\ev')$ such that \textsc{Resample}$(z,\ev')$ was called at setp $t-1$ and $y\in \obj(\ev')\setminus S_z(\ev')$ became bad as soon as the variables
 $\psi_{\obj(\ev')\setminus S_z(\ev')}$  were resampled. But   \textsc{Resample}$(z,\ev')$ must end at a step $t'>t$ and at this step all variables of $\obj(\ev')\setminus S_z(\ev')$ must be good and thus
$y$, which belongs to $\obj(\ev')\setminus S_z(\ev')$, is good at step $t'>t$ in contradiction with the assumption.

%Namely let us suppose that  there exists an event $\ev'$ such that $y\in \obj(\ev')$ which is occuring
%%, $\ev'$ does not occur in the beginning of \textsc{Resample}$(x,\ev)$, it has some of its variables resampled during the execution of \textsc{Resample}$(x,\ev)$
%%and then in
%at the end of \textsc{Resample}$(x,\ev)$. This event $\ev$  may
%have occurred and stopped occurring  in principle several times during the steps of \textsc{Resample}$(x,\ev)$, but if in the end it occurs,
% then there exists a last step $t\le n$ of \textsc{Resample}$(x,\ev)$  such that $\ev'$ was not occurring at step $t-1$, started to occur at step $t$ and kept occurring during the remaining $n-t$ steps of \textsc{Resample}$(x,\ev)$.
%
%\\Consider this last step $t$ of \textsc{Resample}$(x,\ev)$  in which the event $\ev'$ occurs and keeps occurring. There exists an event $\ev''$ and an atom $z$ such that the process \textsc{Resample}$(z,\ev'')$ was called at step $t-1$ of \textsc{Resample}$(x,\ev)$,
% $\{\obj(\ev'')\setminus S_z(\ev''))\}\cap \obj(\ev')\neq\0$ and $\ev'$ occurred as soon as the variables $\psi_{\obj(\ev'')\setminus S_z(\ev''))}$ were resampled. However, the algorithm says that the atoms of  $\{\obj(\ev'')\setminus S_z(\ev''))\}\cap \obj(\ev')$ must be good at the end of \textsc{Resample}$(z,\ev)$, which occurs at some step $t'>t$, and thus at this step $t'$ also $\ev'$ which contains $\varphi$ must not occur and this is in contradiction with the definition of $t$.

~~~~~~~~~~~~~~~~~~~~~~~~~~~~~~~~~~~~~~~~~~~~~~~~~~~~~~~~~~~~~~~~~~~~~~~~~~~~~~~~~~~~~~~~~~~~~~~~~~~~~~~~~~~~~~~~~~~~~~$\Box$

\vskip.2cm

\begin{lem}\label{mtrees}
\textsc{Forest-algorithm} performs at most $m=|\L|$ phases.
\end{lem}

\\{\it Proof}. Consider two phases $l$ and $s$, with $l<s$, generated by an execution of \textsc{Forest-Algorithm} and let $(x_l, \ev_l)$ and $(x_s,\ev_s)$ be the pairs resampled at their initial steps respectively, i.e, the roots of phase $l$ and $s$ respectively.  By Lemma \ref{Agood2}, all atoms  in $\obj(\ev_l)\setminus S_{x_l}(\ev_l)$ are good when phase $l$  ends and at the beginning of any successive phase. In particular, since  $x_l\in \{ \obj(\ev_l)\setminus S_{x_1}(\ev_l)\}$, $x_l$ is good and thus $x_l\notin \obj(\ev_s)$. In conclusion  $x_l \neq x_s$.

~~~~~~~~~~~~~~~~~~~~~~~~~~~~~~~~~~~~~~~~~~~~~~~~~~~~~~~~~~~~~~~~~~~~~~~~~~~~~~~~~~~~~~~~~~~~~~~~~~~~~~~~~~~~~~~~~~~~~~$\Box$

%\\It is important to note that Lemma \ref{mtrees} works since the variable label is never in the seed, so the variable label can be resampled and became good.

\subsection{Witness Forest}\label{swf}

\\We will associate to an execution of \textsc{Forest-algorithm}  a labeled forest  formed by plane rooted trees  whose vertices are labeled with pairs $(x, \ev)$ belonging to $\cal L$.
%We call this forest the {\it witness forest} of the algorithm.

%\begin{defi}[Weakly feasible forest]\label{wff}
%A labeled forest is called weakly feasible if:

%\begin{enumerate}
 % \item  $(v_1, A_1)$ and $(v_2, A_2)$ are labels of roots of two different trees, then $A_1$ and $A_2$ can at most share variables in $S_{v_1^1}(A_1)$.
  %\item $(v_1, A_1)$ and $(v_2, A_2)$ are labels of two siblings, then $A_1$ and $A_2$ can at most share variables in $S_{v_1^1}(A_1)$.
  %\item $(v, A)$ is the label of a node u, then the variable-label of any child of u belongs to $vbl(A) \setminus S_{v}(A)$.
%\end{enumerate}

%\end{defi}

%Indeed, let us justify the name of the algorithm and show that it produces in a natural way a weakly feasible forest.

\\Suppose that the algorithm performs $r$ phases and during  the phase $s$, $s \in \{1, \dots, r\}$, the algorithm performs $n_s$ steps, in such a way that the record of the algorithm is
\begin{equation}\label{forest}
  {\cal L} = \left((x^1_1,\ev_1^1), \dots, (x^1_{n_1},\ev_{n_1}^1),(x^2_1,\ev_1^2), \dots, (x^2_{n_2},\ev_{n_2}^2), \dots, (x^r_1,\ev_1^r), \dots, (x^r_{n_r},\ev_{n_r}^r)\right).
\end{equation}
At each phase $s$, $1\le s \le r$, we will associate a tree $\t'_s$. Let

\begin{equation}\label{t1}
  (x^s_1,\ev_1^s), \dots,({x^{s}_{i}},\ev^s_i), \cdots, (x^s_{n_s},\ev_{n_s}^s),
\end{equation}
be the pairs resampled at phase $s$. We construct the tree $\t'_s$ in the following way.
\vskip.2cm
\\a) The root of  $\t'_s$ has label $(x^s_1,\ev_1^s)$.

\\b) For $i>1$, we proceed  by checking  if $({x^{s}_{i}},\ev^s_i)$ is such that ${x^{s}_{i}} \in (\obj(\ev^{s}_{i-1}) \setminus S_{x^{s}_{i-1}}(\ev^s_{i-1}))$,

\\- if yes, we add $({x^{s}_{i}}, \ev^s_i)$ as a child of $({x^{s}_{i-1}},\ev^s_{i-1})$,

\\- if no, we go back in (\ref{t1}) checking the {\it ancestors} of the vertex labeled by $({x^{s}_{i-1}},\ev^s_{i-1})$  until we find a pair $({x^{s}_{j}},\ev^s_j)$, with $j<i$, such that ${x^{s}_{i}} \in (\obj(\ev^s_{j}) \setminus S_{x^{s}_{j}}(\ev^s_{j}))$, and we add $({x^{s}_{i}},\ev^s_i)$ as a child of $({x^{s}_{j}},\ev^s_{j})$.

%\begin{rema}
 % Let $v$ and $u$ be vertices in the rooted tree $\tau$ such that the distance of the root to $v$ is $i$ while the distance of the root to $u$ is $j$, with $j>i$. We say that the vertex $v$ is an antecedent of $u$ if there is a path from $v$ to $u$ with exactly $j-i$ edges.
%\end{rema}

\\Observe that by the construction of \textsc{Forest-Algorithm}, all pairs $({x^{s}_{i}},\ev^s_i)$ can be added to $\t'_s$ in this way, then $\t'_s$ has $n_s$ vertices (leaves included) with labels $({x^{s}_{i}},\ev^s_i)$ with $i=1,\dots, n_s$. By Lemma \ref{Agood2} the pair $(x^{s+1}_1,\ev_1^{s+1})$ is the first pair in (\ref{forest}) that can not be added to $\t'_s$ in this way, so we build a new tree $\t'_{s+1}$ with root $(x^{s+1}_1,\ev_1^{s+1})$ following the same rule described to build $\t'_s$.

%This definition of the witness forest of the algorithm together with Theorem \ref{Agood} immediately implies the following corollary.
\\Note that the vertices of the forest defined above are naturally ordered according to the natural order of the steps made by  the algorithm. The  forest ${ F}'=\{\t'_1,\dots, \t'_r\}$ so obtained uniquely associated to the record $\mathcal{L}$ is such  that, for each $s\in [r]$,   $\t'_s$ is a rooted plane tree with  $n_s$  vertices and each vertex of $\t'_s$ has label $(x,\ev)$ where $x \in \obj(\ev)$ and $\ev \in \FF$.

\\Note that, by Lemma \ref{mtrees} we have that $r\le m$ and thus the   forest ${ F}'$ contains at most $m$ trees.

\\Note also that in each tree $\t'_s$ of $F'$ the list of labels of the vertices of $\t'_s$  ordered according to the depth-first search, coincides with  the list  (\ref{t1}).

\\Note finally that,
by construction, the correspondence ${\cal L}\mapsto F'$ is an injection.

%\\We now  prove two lemmas that show some important properties on the forest generated by the algorithm.

\begin{lem}\label{dr}
Consider a tree $\t' \in F'$, and let $v_i$ and $v_j$ be two vertices in $\t'$ with labels $(x_i, \ev_i)$ and $(x_j, \ev_j)$ respectively. We have that

\begin{enumerate}

\item[\rm a)] If $v_i$ is a child of $v_j$, then $x_i \in \obj(\ev_j) \setminus S_{x_j}(\ev_j)$.

\item[\rm b)] If $v_i$ and $v_j$ are siblings in $\t'$, then $x_i\neq x_j$.% and $ \left(vbl(A_i)\cap vbl(A_j)\right)\subseteq S_\psi(A_i)$.
\item[\rm c)]  Any vertex $v \in \t'$ with label $(x, \ev)$ has at most $\|\ev\|$ children, where $\|\ev\|$ is defined in (\ref{power2}).
\end{enumerate}
\end{lem}

\\{\it Proof.}

\\a) It is trivial by construction of the algorithm.
\vskip.2cm
\\b) As $v_i$ and $v_j$ are siblings, suppose that $v_i$ and $v_j$ are the i-th and the j-th children of a vertex in $\t'$, with  $i < j$ in the natural order of the vertices of $\t'$ induced by the steps of the algorithm. For $q$ such that $i\le q<j$, let $(x_q, \ev_q)$ be the label of the $q^{th}$ sibling. By Lemma \ref{Agood2} when \textsc{Resample}$(x_q,\ev_q)$ ends all the atoms $x_i,\dots, x_q$ are good as well the atoms in  $(\obj(\ev_i)\setminus S_{x_i}(\ev_i))\cup\dots\cup (\obj(\ev_q)\setminus S_{x_q}(\ev_q))$. Therefore $x_j$ can not be in the set $\{x_i, x_{i+1},\dots, x_{j-1}\}$.
\vskip.2cm
\\c) Follows trivially from items a) and b).

~~~~~~~~~~~~~~~~~~~~~~~~~~~~~~~~~~~~~~~~~~~~~~~~~~~~~~~~~~~~~~~~~~~~~~~~~~~~~~~~~~~~~~~~~~~~~~~~~~~~~~~~~~~~~~~~~~~~~~~~$\Box$

\vskip.2cm

\\Given a forest $F'$ produced by the algorithm, we let $X_{F'}$ be the set of atoms which label the roots of the trees of the forest, i.e.,
$$X_{F'}=\{ x \in \L:~ \exists  \ev \in \FF \ \mbox{such that} \ (x, \ev)\ \mbox{is the root label of some}\; \tau' \in F' \}.$$
Lemma \ref{mtrees} implies  that atoms  in $X_{F'}$ are all distinct.

\begin{defi}[Witness forest]
Given the record $\mathcal{L}$ of \textsc{Forest-Algorithm} and the forest $F'$ associate to $\mathcal{L}$, we construct a new forest $F$ by adding to $F'$ new vertices in the following way:

  \\1) Add to the forest $F'$ as many isolated vertices  as the atoms which are in $\L\setminus {X}_{F'}$, and  give to these isolated vertices the label $(x, \emptyset)$ for all $x \in \L\setminus {X}_{F'}$.

  \\2) For each vertex $v$ of the forest $F'$ with label $(x, \ev)$ with less than $\|\ev\|$ children, do the following: let $H_v$ be the set of atoms  in $\obj(\ev)\setminus S_x(\ev)$ which are not atoms labels of the children of $v$. For each $y \in H_v$ we add to $v$ a leaf with label $(y,\0)$ in such a way that $v$ has now exactly  $\|\ev\|$ children.

\end{defi}

\\The  new labeled forest $F$, so obtained uniquely  associated to the random variable $\mathcal{L}$ by the prescriptions described above, is called the {\it witness forest} produced by \textsc{Forest-Algorithm}. This witness forest $F$   has, by construction,  the following properties.
\vv

\\{\it  Properties of the witness forest $F$.}

  \begin{enumerate}

\item $F$  is constituted by exactly $|\L|=m$ labeled rooted trees $\t_1,\dots,\t_m$ (some of which are just isolated vertices).

\item Let the vertex $u$ be a child of the vertex $v$ in $\t \in F$ and let $(x_u, \ev_u)$ and $(x_v, \ev_v)$ be their labels respectively. Then $x_u \in \obj(\ev_v)$.

\item  Each internal vertex $v$ of $\t \in F$ carries a label $(x_v, \ev_v)$ where $x_v\in \obj(\ev_v)$ and $\ev_v\in \FF$, while each leaf $\ell$ of $\t$ carries a label $(x_\ell,\0)$ and $x_\ell\in \obj(\ev_w)$, where $w$ is the vertex parent of $\ell$.

\item Let the vertices $v$ and $v'$ be the $i^{th}$ and the $j^{th}$  siblings in $\t \in F$, with $i<j$ in the depth-first search order of $\t$, and let $(x_i, \ev_i)$ and $(x_j, \ev_j)$ be their labels respectively, then $x_i\neq x_j$.

\item Let $v$ be an internal vertex of $\t \in F$  and let $(x, \ev)$ be its label, then the vertex $v$ has exactly $\|\ev\|$ children.
\end{enumerate}

%\begin{rema}
%Thus, the witness forest of the algorithm  is such that $F^*$  has exactly $m$ trees, some of which are just isolated vertices, and any  internal vertex  $v$ of $F^*$, with event label $A_v$, has $|A_v|-\|A_v\|$ children and there is a one-to-one correspondence between internal vertices of $F^*$ and vertices of $F$ (i.e. all the internal vertices of $F^*$  are vertices of $F$, while all the isolated roots or leaves are not vertices of the forest $F$).
%\end{rema}

\\Let  $\mathcal{F}_n$ be the set of labeled  forests satisfying properties 1-5 above, that contains  $n$ internal vertices in total and let
$\mathcal{F}=\cup_{n\ge 0} \mathcal{F}_n$.

\\It is important to stress that the map ${\cal L}\mapsto F$ is an injection. Therefore, since \textsc{Forest-Algorithm} lasts $n$ steps if and only if the witness forest associated to the record $\cal L$ of \textsc{Forest-Algorithm} has $n$ internal vertices, the probability $P_n$ defined in (\ref{p1}) can be written as
\begin{equation}\label{wf}
P_n    = \prob(\mbox{the witness forest associated to $\cal L$ has} \ n \ \mbox{internal vertices}).
\end{equation}

\\Then the next goal is to estimate the probability that \textsc{Forest-Algorithm} produces a witness forest $F$ with $n$ internal vertices.

\subsection{The validation algorithm}\label{ssc}

\begin{defi}[Admissible sequence]
We say that a sequence $S=\{(x_1,\ev_1), \cdots, (x_n, \ev_n)\}$ is admissible if $x_i \in \obj(\ev_i)$, for all $i=1, \cdots, n$.
\end{defi}

\\Given a witness forest $F$ with $n$ internal vertices, we can associate to $F$, in a natural way, the  admissible sequence $S_F=\{(x_1,\ev_1), \dots, (x_n,\ev_{n})\}$  formed by the labels of its internal vertices. Namely, the sequence $S_F$ coincides with (\ref{forest}).

\\We now describe a validation algorithm,  called \textsc{S-Check}, whose input is an admissible sequence $S=\{(x_1,\ev_1), \dots, (x_n,\ev_{n})\}$. \textsc{S-Check} first samples all variables in $\psi_\L$ and then resamples some of the variables in $\psi_\L$.

\vskip.5cm
\fbox{
\begin{minipage}{0.9\linewidth}
{\textsc{ S-Check}}.
\vskip.2cm
Given the admissible sequence $S=\{(x_1,\ev_1), \dots, (x_n,\ev_{n})\}$
\begin{enumerate}
  \item Sample all variables in $\psi_\L$.
  \item  For $i= 1, \cdots, n$, do
  \item   If $\ev_{i}$ occurs, resample all the variables $\psi_y$ with $y \in \obj(\ev_{i})\setminus S_{x_i}(\ev_i)$. If the event $\ev_{i}$ does not occur, return {\it failure}.
  \item End for.
  \end{enumerate}
\end{minipage}
}
\vskip.5cm
\\The procedure described at line 3 of \textsc{ S-Check} is called a {\it step}. Of course, if $S=\{(x_1,\ev_1), \dots, $ $(x_n,\ev_{n})\}$ is the input for \textsc{ S-Check},  its execution will perform exactly $n$ steps if it does not return failure. Observe that  \textsc{S-Check} does not return failure (i.e. passes) if, and only if, in each step $i$ the event $\ev_i$ occurs under the current evaluation of the variables.
%\vskip.2cm
%\\{\bf Remark}.
%It is important to stress that in the $i$-th repetition of line 3 of \textsc{ S-Check} the only information that we have is if the event $\ev_i$ occurs or not, in fact we do not know the evaluation of the variables in $\obj(\ev_i)$.
%\vskip.2cm

\begin{lem}\label{Fcheck1}
Let   $S=\{(x_1,\ev_1), \dots, (x_n,\ev_{n})\}$ be an admissible sequence. Then
\begin{equation}
\prob(\mbox{\textsc{$S$-check} with input $S$ passes})\leq \prod_{i=1}^{n} \prob(\ev_i).
\end{equation}
\end{lem}

\\{\it Proof}. Consider the first step of \textsc{$S$-check}: we sample all variables $\psi_\L$ reaching a configuration $\o_0$ and we have to check if the event $\ev_1$ happens, if $\ev_1$  does not happen we stop, otherwise we resample the variables in $\obj(\ev_1) \setminus S_{x_1}(\ev_1)$. As $S_{x_1}(\ev_1)$ is a seed, we have that the new configuration $\omega_1$ is such that any event $\ev$ has probability to occur at most $\prob(\ev)$. Therefore, by induction, at each step $i$ the probability of $\ev_i$ to occur is at most $\prob(\ev_i)$. As \textsc{$S$-Check} is successful if and only if all events $\ev_i$ occur, then
\begin{equation}
\prob(\mbox{\textsc{$S$-check} with input $S$ passes})\leq \prod_{i=1}^{n} \prob(\ev_i).
\end{equation}

\begin{lem}\label{Fcheck}
Given a witness forest $F\in \mathcal{F}_n$ whose  internal vertices carry labels
$$
S_F=\{(x_1,\ev_1), \dots, (x_n,\ev_{n})\},
$$
we have that
\begin{equation}\label{pfsc}
\prob(\mbox{\textsc{Forest-Algorithm} produces} \ F) \leq \prod_{i=1}^{n} \prob(\ev_i).
\end{equation}
\end{lem}

\\{\it Proof.}
Observe that if all the random choices made by an execution of \textsc{Forest-Algorithm} that produces $F$ as witness forest are also made by the Algorithm \textsc{$S$-Check} with input $S_F$, then in each step $i$ the event $\ev_i$ occurs and so \textsc{$S$-Check} does not return failure. Then
\begin{equation}
\prob(\mbox{the \textsc{Forest-Algorithm} produces  $F$}) \leq \prob(\mbox{\textsc{$S_F$-Check} with input $S_F$ passes}).
\end{equation}
Now (\ref{pfsc}) follows from Lemma \ref{Fcheck1}

~~~~~~~~~~~~~~~~~~~~~~~~~~~~~~~~~~~~~~~~~~~~~~~~~~~~~~~~~~~~~~~~~~~~~~~~~~~~~~~~~~~~~~~~~~~~~~~~~~~~~~~~~~~~~~~~~~~~~~~~$\Box$

%\begin{lem}[Randomness Lemma]\label{rl}
%At the beginning of any repetition of the line 2 of $F$-check, the distribution of the current assignment of values of the variables in $\Psi$ is as if all variables have been sampled anew.
%\end{lem}

%{\bf Proof:}
%{\red see Harvey Lemma 19. }
%It is important to note that in the phase $i$ of the F-check just the values of the variables in $vbl(A_{i})\setminus S_{A_i}$ are resampled, then, at the beginning of phase $i+1$ the distribution of the evaluation of the variables is the same as the Line 1. Indeed, the fact that the variables in $S_{A_i}$ are not resampled does not change the probability of $A_i$ occurring and neither the probability of any other event occurring, as by the Definition \ref{seed} the knowledge that the event $A_i$ happens says nothing about its seed, we can have any configuration for the seed.$\Box$

% (Randomness Lemma). At the beginning of any iteration of the for loop of Step 2 of ValAlg, the distribution of the current assignment of values to the variables Xi, i = 1, ..., l is as if all variables have been sampled anew. Therefore the probability of any event occurring at such an instant is bounded from above by p.

\vskip.2cm

\\{\bf Remark}.
In the $i$-th phase of the \textsc{$S$-Check} just the values of the variables $\psi_{\obj(\ev_{i})\setminus S_{x_i}(\ev_i)}$ are resampled.  So, at the beginning of phase $i+1$ the distribution of the evaluation of the variables is the same as the Line 1. This is not the case for \textsc{Forest-Algorithm}, once it would mean that this algorithm is not making any progress in the search of a configuration such that any event in $\FF$ occurs.

\subsection{The unlabeled Forest}

\\The strategy to prove Theorem \ref{teo1} is to show that the probability that \textsc{Forest-Algorithm} lasts at least $n$ steps decays exponentially in $n$, which implies that \textsc{Forest-Algorithm} terminates almost surely, returning an evaluation of $\psi_\L$ such that all events in $\FF$ do not occur.

\\If \textsc{Forest-Algorithm} lasts $n$ steps  then it produces a witness forest  with $n$ internal nodes. Recall  that, if the internal vertex $v$  of the witness forest has event label $\ev_v$  then this vertex has exactly $s_v = \|\ev_v\|$ children.

\\Let ${\cal F}^*_n$ be  the set of all unlabeled forests constituted by $|\L|=m$ plane trees  having in total $n$ internal vertices and such that each internal  vertex $v$ has a number of children in the set   $E'_{\FF}$ defined in (\ref{E'}).
Given the record $\cal L$ of \textsc{Forest-Algorithm} such that $|{\cal L}|=n$ and given the witness forest $F\in \mathcal{F}_n$ associated to $\cal L$, we define the function
$$f: \mathcal{F}_n \rightarrow {\cal F}^*_n$$
 $$ f(F)=\Phi$$
that removes all the labels of $F$ obtaining an unlabeled witness forest $\Phi \in {\cal F}^*_n$. We call $\Phi$ the {\it associated unlabeled witness forest} produced by \textsc{Forest-Algorithm}.  Given an internal vertex $v$ of an unlabeled forest $\F\in {\cal F}^*_n$, we let 
 $s_v$ be the number of children on $v$.

\\For $\F \in {\cal F}^*_n$ let us define
\begin{equation}
  P_\F = \prob( \F \ \mbox{is the associated unlabeled witness forest produced by \textsc{Forest-Algorithm}}).
\end{equation}

%\\Given the witness forest $F$, let $g(F)$ be the associated unlabeled forest. Note that if the internal vertex $v$  of $F$ has event label $A$ such that $|A|=k$ and $\|A\|=l$, then this vertex has exactly $k-l$ children. Thus we can also construct a labeled forest $f(F)$ having the same underlying unlabeled forest of $F$, i.e. $g(F)$, in such a way that if $v$ is an internal vertex of $F$ with event label  $A$ with $|A|=k$ and $\|A\|=l$, then the same vertex in $f(F)$ has label $(k,l)$ and has $k-l$ children. Let us denote the forest $f(F)$ by {\it weakly labeled forest}.

%\\Let thus $\F$ be a weakly labeled  forest with $n$ internal vertices, each vertex $v$ with labels $(k_v,l_v)$ and thus with $k_v-l_v$ children, and let us define
%\begin{equation}
%P_\F= \mathbb{P}(\F \ \mbox{is the weakly labeled forest associated to the witness forest}).
%\end{equation}

%\\Let ${\cal F}_n$ be the set of the labeled forests with $n$ internal vertices, we have that
\begin{eqnarray}
  P_\F &=& \sum_{F\in {\mathcal{ F}}_n ; \atop f(F)=\F}\prob(\mbox{the \textsc{Forest-Algorithm} produces the witness forest} \ F) \label{1}\\
       &\le& \sum_{F\in {\mathcal{F}}_n ;\atop f(F)=\F}\prod_{v\in F}\prob(\ev_v) \label{2}\\
       &\leq& \prod_{v\in \F}p_{s_v}\sum_{F\in {\mathcal{F}_n};\atop f(F)=\F}1, \label{3}
\end{eqnarray}
where $p_s$ is defined in (\ref{pij}) and inequality (\ref{2}) is due to Lemma \ref{Fcheck}. Now observe that
\begin{equation}
\sum_{F\in {\mathcal{ F}_n;}\atop f(F)=\F}1\le \prod_{v\in \F}d_{s_v},
\end{equation}
since for each vertex $v \in \F$ with $s_v$ children, we have $d_{s_v}$ options for its event label, and fixed the atom label and the event label of the parent of $v$, we can determine uniquely the atom label of $v$. Indeed, suppose that the parent of $v$ is the vertex $u$, which has $s_u$ children, and $x_u$ is its atom label, then for the event label of $u$ we have $d_{s_u}$ options, and once fixed the event label of $u$, suppose $\ev_u$, we know the atom labels of all children of $u$, namely they are, in order, the atoms in $\obj(\ev_u) \setminus S_{x_u}(\ev_u)$. So, if $v$ is the $i$-th child of $u$, then its atom label is the $i$-th atom in $\obj(\ev_u) \setminus S_{x_u}(\ev_u)$. Proceeding recursively, observe that what we need to know is the atom label of the roots of each tree in $\F$, however this information is easily obtained by the construction of a witness forest, as the trees are organized by the atom labels of their roots.

\\Then,
\begin{equation}
P_\F \leq \prod_{v\in \F}d_{s_v}p_{s_v}.
\end{equation}

\\We now can bound  the probability $P_n$ (see (\ref{wf})) that \textsc{Forest-Algorithm} lasts $n$ steps as
$$
P_n\le  \sum_{\F\in {\cal F}^*_n}{P_\F}.
$$
\\To estimate  $\sum_{\F\in {\cal F}^*_n}{P_\F}$, observe that every forest $\F\in {\cal F}^*_n$ is constituted by $m$ trees $\t_1,\dots,\t_m$
with $n_1,\dots, n_m$ internal vertices respectively. The numbers $n_1, \dots,n_m$ are such that $n_i\ge 0$ for all $i=1,\dots, m$ and $n_1+n_2+\dots+n_m=n$.
Recall also  that the number of children of the internal vertices  of any  $\t_i$ takes values in the set $E'_\FF$.
Let us denote by $\cal T$ the set of plane trees with number of children of the internal vertices taking  values in the set $E'_\FF$ and let ${\cal T}_n$ be the subset of  $\cal T$ formed by the trees with exactly $n$ internal vertices.

\\Let us denote  shortly,  for $s \in E'_\FF$,
\be\label{wk}
w_s= d_{s}p_{s}.
\ee
For  a tree $\t\in \cal T$, let $V_{\t}$ be the set of its internal vertices. Then  define the weight of $\t$ as
$$
\o(\t)= \prod_{v\in V_\t}w_{s_v}
$$
where we recall that $s_v$ is the number of children of the vertex $v$.  

\\For a given $n\in \mathbb{N}$, let
$$
Q_n=\sum_{\t\in {\cal T}_n}\o(\t).
$$
Therefore, the probability that the \textsc{Forest-Algorithm} lasts $n$ steps is bounded by
\be\label{penne}
P_n\le  \sum_{n_1+\dots +n_m=n\atop n_i\ge 0}Q_{n_1}\dots Q_{n_m}.
\ee
It is now easy to check that  $Q_n$ is defined by the recurrence relation

\be\label{recur1}
Q_n=\sum_{s\in E'_{\FF}}w_s\sum_{n_{1}+\dots +n_{s}=n-1\atop n_1\ge 0, \dots ,n_{s}\ge 0}
Q_{n_1}\dots Q_{n_{k-l}},
\ee
with $Q_0=1$.
Now let
$$
W(z)=\sum_{n=1}^\infty Q_n z^n,
$$
be the generating function encoding the sequence $\{Q_n\}_{n\ge 1}$. Then
we have  from (\ref{recur1})
\begin{eqnarray}
  W(z) &=&  z\sum_{n=1}^{\infty}  \sum_{s\in E'_{\FF}}w_s\sum_{n_{1}+\dots +n_{s}=n-1\atop n_1\ge 0, \dots ,n_{s}\ge 0}
Q_{n_1}z^{n_1}\dots Q_{n_{s}}z^{n_{s}}  \\
       &=&  z \sum_{s\in E'_{\FF}}w_s\prod_{i=1}^s\sum_{n_i\ge 0}Q_{n_i}z^{n_i} \\
       &=& z \sum_{s\in E'_{\FF}}w_s\prod_{i=1}^s[1+ \sum_{n_i\ge 1}Q_{n_i}z^{n_i}] \\
       &=& z  \sum_{s\in E'_{\FF}}w_s(1+W(z))^s,
\end{eqnarray}
i.e. denoting, for $\x>0$
\begin{equation}
\phi_{\FF}(\x)=\sum_{s\in E'_\FF}w_s(1+\x)^s %= \sum_{k\in \EA}\sum_{l\in L_k} d_k^l p_k^l (1+\x)^{k-l},
\end{equation}
we have
\begin{equation}
W(z)= z \phi_\FF(W(z))
\end{equation}
and thus, by a well known result in analytic combinatorics (see e.g. Proposition IV.5 of \cite{FS} or also Theorem 5 in  \cite{Dm}) we have that
the coefficients of the generatin function $W(z)$ are bounded as follows.
\begin{equation}
Q_n\le \r^n
\end{equation}
where
\begin{equation}
\r=\min_{\x>0}{\phi_\FF(\x)\over \x}.
\end{equation}

\\Hence,
\begin{equation}
P_n\le \sum_{n_1+\dots +n_m=n\atop n_i\ge 0}Q_{n_1}\dots Q_{n_m}\le \r^n  \sum_{n_1+\dots +n_m=n\atop n_i\ge 0}1=\r^n {n+m-1\choose m-1}
\end{equation}
Now,  if condition (\ref{teo}) holds, we have that the probability that the \textsc{Forest-Algorithm} runs at least  $n$ steps decays exponentially in $n$ if $n$ is sufficiently large.  In particular it is easy to check that
$$
\r^n {n+m-1\choose m-1}\le \r^{n\over 2}
$$
as soon as
$$
{n\over \ln n}\ge {2m\over |\ln(\r)|},
$$
i.e.
as soon as
$$
n\ge {2m\over |\ln(\r)|}\ln^2\Big({2m\over |\ln(\r)|}\Big) \equiv N.
$$
Thus, if we estimate $P_n=1$ if $n\le N$ and $P_n\le \r^{n/2}$ if $n>N$, the expected number of steps $T$ of \textsc{Forest-Algorithm} is given by
$$
T\le {N(N+1)\over 2}+ \sum_{n=N+1}^\infty n\r^{n\over 2}.
$$

\section{Examples}\label{appl}

\\In what follows $G=(V,E)$ is  a graph  with maximum degree $\D$ and   $k \in \mathbb{N}$. A coloring of the vertices (resp. edges) of $G$ is a function $c: V \to [k]$ (resp. $c': E \to [k]$).

\vskip.2cm
\\\underline{{\it Example 1: Nonrepetitive vertex coloring of a graph}}
\vskip.2cm
\\A coloring of the vertices of $G$ is  {\it nonrepetitive} if, for any $n \ge 1$, no path $p=\{v_1,v_2,\dots, v_{2n}\}$ is colored repetitively, i.e. such that $c(v_i)=c(v_{i+n})$ for all $i=1,2,\dots n$. The minimum number of colors needed such that $G$ has a non-repetitive vertex coloring
 is called the non repetitive chromatic index of $G$ and it is denoted by $\p(G)$. Here we  are in the uniform variable setting, where
the set of atoms  $\L$ coincides with $V$ and to each atom/vertex $v\in V$ we associate a random variable
$\psi_v$, the color of $v$, that takes values in $[k]$ according to the uniform distribution.
Let $P_{n}$ be the set of all paths with $2n$ vertices and set $P= \cup_{n \geq 1}{P_n}$. The family $\FF$ of bad events is the set
$\FF=\{\ev_{p}\}_{p \in P}$, where $\ev_p$ is the event ``the path $p$ is colored repetitively''.
For any $p\in P$, given a vertex $v$  a seed of $\ev_p$ not containing $v$ is the half of $p$ that does not contain $v$. Thus, 
if $p$ is a path with $2n$ vertices,  then $\ev_p$ is tidy with seeds of size $n$ and therefore  $\|\ev_p\|=n$  and  $\prob(\ev_p)={1\over k^n}$. So in this case
$$
E'_\FF=\{1,2,3,\dots\}.
$$
In order to apply Theorem \ref{teo1} we have to estimate $d_s$, the maximum number of events of power $s$  containing a fixed vertex. In the present case $d_s$ coincides with the maximum number of paths in $G$ of size $2s$ containing a fixed vertex.  We have 
$$
d_s \leq s\D^{2s
-1}.
$$
Therefore the function $\phi_\FF(\x) $ defined in (\ref{fea}) is in the present case
\begin{eqnarray}
\nonumber  \phi_\FF(\x) &=& \sum_{s\ge 1} s\D^{2s-1} {1\over k^s}(\x+1)^{s} \\
\nonumber              &=& {1\over \D}\sum_{s\ge 1} s\left({\D^2\over k}(\x+1)\right)^s\\
\nonumber              &=& {1\over \D} {{\D^2\over k}(\x+1)\over \left(1-{\D^2\over k}(\x+1)\right)^2}\\
\nonumber              &=& {1\over \D}{(b+1)(\x+1)\over (b-\x)^2},
\end{eqnarray}
where in the last line we have set
$$
k=(1+b)\D^2.
$$
Thus condition (\ref{teo}) is in this case
$$
\min_{\x>0} \left({1\over \D}{(b+1)(\x+1)\over \x(b-\x)^2}\right)<1.
$$
Observe that the minimum occurs at
$$
\x_0={\sqrt{9+8b}-3\over 4},
$$
and
$$
{\phi_\FF(\x_0)\over \x_0}={1\over \D}{\sqrt{(8b+9)^3}+8b^2+ 36b+27\over 8b^3},
$$
and thus if we let $b_0(\D)$ be the solution of the equation
$$
{\sqrt{(8b+9)^3}+8b^2+ 36b+27\over 8b^3}=\D,
$$
we have that the non repetitive chromatic index $\pi(G)$ of a graph with maximum degree $\D$ is such that
\be\label{a}
\pi(G)\le (1+b_0(\D))\D^2.
\ee

%{\red tirar}
%Let us compare with Theorem 8 in \cite{DJKW}which states
%$$
%\pi(G)\le \Big(1+c_0(\D))\Big)\D^2
%$$
%with
%$$
%c_0(\D)={1\over \D^{1\over 3}-1}+ {1\over \D^{1\over 3}}
%$$
%We have e.g. for $\D=3$
%$$
%b_0(3)=2.3143~~~~~~~~~~~~~~~~~c_0(3)=2.9545
%$$
%$$
%b_0(10)=1.29~~~~~~~~~~~~~~~~~c_0(10)=1.33
%$$
%$$
%b_0(100)=0.487~~~~~~~~~~~~~~~~~c_0(100)=0.49
%$$
%$$
%b_0(10^3)=0.205~~~~~~~~~~~~~~~~~c_0(10^3)=0.21
%$$
%$$
%b_0(10^6)=0.019~~~~~~~~~~~~~~~~~c_0(10^6)=0.0201
%$$
%and
%$$
%\lim_{\D\to \infty} \D^{1/3}b_0(\D)=1.88988~~~~~~~~\lim_{\D\to \infty} \D^{1/3}~c_0(\D)=2
%$$
%So we see that our bound is slightly better than \cite{DJKW}.

\\Comparing our bound with Theorem 8 in \cite{GMP},  which states that
\be\label{b}
\pi(G)\le \D^2+\D^{3\over 2}\left[{3\over 2^{2/3}}+ {2^{2/3}\over \D^{1\over 3}-2^{1\over 3}}\right]
\ee
%with
%$$
%d_0(\D)=\left[{3\over 2^{2/3}}+ {2^{2/3}\over \D^{1\over 3}-2^{1\over 3}}\right] \D^{-{1\over 3}}
%$$
%We have
%$$
%b_0(3)=2.3143~~~~~~~~~~~~~~~~~d_0(3)=7.3469
%$$
%$$
%b_0(10)=1.29~~~~~~~~~~~~~~~~~d_0(10)=1.701
%$$
%$$
%b_0(100)=0.487~~~~~~~~~~~~~~~~~d_0(100)=0.5082
%$$
%$$
%b_0(10^3)=0.205~~~~~~~~~~~~~~~~~d_0(10^3)=0.207
%$$
%$$
%b_0(10^6)=0.019~~~~~~~~~~~~~~~~~c_0(10^6)=0.019
%$$
%and
%$$
%\lim_{\D\to \infty} \D^{1/3}b_0(\D)=1.88988~~~~~~~~\lim_{\D\to \infty} \D^{1/3}~d_0(\D)=1.88988
%$$
we observe that  bound (\ref{a}) is better than (\ref{b}) for low values of $\D$ while becomes asymptotically equivalent for
large values of $\D$.

  \vskip.5cm
\\\underline{\it Example 2: Facial Thue Choice Index of planar graphs}
\vskip.2cm
\\We suppose here that the graph $G=(V,E)$ is planar. Suppose moreover that for all edge $e \in E$, a  list $L_e$ of $k$ colors is given. A facial path of  $G$ is a path of $G$ which is part of the  boundary of a face of $G$. The least integer $k$ such that for every collection of lists  $\{L_e\}_{e\in E}$ with $|L_e|=k$ there is an edge coloring of $G$ such that every facial path of $G$ is nonrepetitive is called the {\it  facial Thue choice index of $G$} and is denoted by $\pi'_{fl}(G)$. Observe that the set of independent random variables is in this case $\Psi=\{L_e\}_{e\in E}$.

\\Let $P$ denotes the set of all facial paths with even number of edges. For all $p\in P$ let $\ev_p$ be the event  ``$p$ is repetitive'', i.e.,  if $p=\{e_1,\dots,e_n,e_{n+1},\dots, e_{2n}\}$ we have $c'(e_i)=c'(e_{i+n})$ for all $i\in [n]$ where $c'(e)$ is the color chosen in the list $L_e$ via the random experiment. The family of bad events is thus
$\FF=\{\ev_p\}_{p\in P}$.
Observe that, analogously to the previous example, any event $\ev_p$ with $p\in \cal P$ is tidy and we can take as a seed of $\ev_p$ the first or the second half of the path $p$. This implies that as before $E'_\FF=\{1,2,3,\dots,\}$. Moreover, if  $|p|=2n$, we have that $\prob(\ev_p)\le {1\over k^n}$, and since every edge of a planar graph is contained in at most $4n$ facial paths of $G$ of size $2n$, for $s\in E'_\FF$ we have that
$ d_s\le 4s,$
and therefore
\begin{eqnarray}
\nonumber  \phi_\FF(\x) &\leq& \sum_{s\ge 1}{1\over k^s}4s(\x+1)^s \\
\nonumber              &<& {4{\x+1\over k}\over \left(1-{\x+1\over k}\right)^2} \\
\nonumber              &=& {4k(\x+1)\over (k-\x-1)^2}.
\end{eqnarray}
Then, we have
$$
\min_{\x>0}{\phi_\FF(\x)\over \x}<1
$$
as soon as $k\ge 12$, which is the same bound obtained in \cite{Pr} via entropy compression method.

  \vskip.5cm
\\\underline{\it Example 3: Coloring graphs frugally}
\vskip.2cm

\\A proper vertex coloring of a graph $G$ is said $\b$-frugal if any vertex has at most $\b$  members of any color class in its neighborhood. The minimum number of colors required such that a graph $G$ has at least one $\b$-frugal proper vertex coloring is called the $\b$-frugal chromatic number of $G$ and will be denoted by $\chi_\b(G)$.
Analogously to the Example 1,  we  are in the entropy compression setting  where
 $(\L,k)\equiv (V,k)$ and to each $v\in V$ we associate a random variable
$\psi_v$  (the color of $v$) that takes values in $[k]$ according to the uniform distribution.

%\vskip.2cm
%\\Let $G=(V,E)$ be a graph with maximum degree  $\D$. Suppose that we have $N$ colors ordered $\{1, 2, \cdots, N\}$. Let the set of the random variables $\Psi=\{c(v)\}_{v\in V}$ be the colors of each vertex. The colors are assigned with the uniform distribution.

\\Observe that in the present case we have only two kind of bad events. First the coloring has to be proper. So,  for each edge $e=\{u,v\}$ of $G$ we must avoid the event $\ev_e$ that  ``$u$ and $v$ have the same color'', and let $\FF_1=\{\ev_e\}_{e \in E}$.

\\We say that a set $\s$ formed by  $\b+1$ vertices of $G$ is a $\b$-star  of $G$ if all members of $\s$ are
 neighbors of a common vetex $v\in V$, in other words, if there is $v\in V$ such that
$\s\subset \G_G(v)$. Let $S_\b$ denote the set of all  $\b$-stars of $G$.
Given $\s\in S_\b$, let $\ev_\s$ be the event ``all the $\b+1$ vertices forming  $\s$ receive the same color'', i.e., $\s$ is monochromatic. We thus have a second family of bad events $\FF_\b= \{\ev_\s\}_{\s\in S_\b}$. Clearly the events of the family $\FF_1$ are tidy with seeds
of size 1 and power equal to 1, while all events of the family $\FF_\b$ are tidy with seeds of size 1 and power equal to $\b$. A $\b$-frugal coloring of the vertices of $G$ occurs if none of the events of the family $\FF=\FF_1\cup \FF_\b$ occurs.

%\\Observe that for every $e \in E$ we have that $|\ev_e|=2$ and $\|A_e\|=1$, since if we know that the edge $e=\{u,v\}$ is monochromatic, letting the color of $u$ or $v$ without resampling we do not help any event in $A$ or $B$ to happen, once we do not know the color of this  vertex. Indeed, by Remark \ref{od}, the number of configurations that makes $A_e$ happens is $N$, while $N^{|S|}=N$, where $S$ is one of the vertices of $e$. Then the seed of $A_e$ is composed by one of the vertices of $e$ and obviously $S_{u}(A_e)=\{v\}$ and $S_{v}(A_e)=\{u\}$.
%
%\vskip.2cm
%\\On the other hand, for every  $v \in V$ and $h_{\b}^{v} \in H_{\b}^{v}$ we have that $|B_{h_\b^v}|= \b+1$  and $
%\|B_{h_\b^v}\|=1$, since letting the color of one vertex in ${h_\b^v}$ without resampling we do not have any information and so do not help any event in $A$ or $B$ to happen, while letting two vertices without resampling we could have an edge monochromatic and then it would help one of the events in $A$ to happen. Indeed,  by Remark \ref{od}, the number of configurations that makes $|B_{h_\b^v}$ happens is $N$, while $N^{|S|}=N$, where $S$ is one of the vertices of $h_\b^v$. Its is clear that we can always choose a seed without a vertex, for example, if $u \in h_\b^v$ we can set that
%\begin{equation}
%S_u(B_{h_\b^v})= \{u'; \mbox{where} \ u' \ \mbox{is the smallest vertex greater than} \ u\}.
%\end{equation}

\\In the present case, $E'_\FF=\{1, \b\}$,
and for  every $e \in E$ and $\s \in S_\b$, we have that $\prob(\ev_e)= \frac{1}{k}$ and $\prob(\s) \leq \frac{1}{k^{\b}}$ respectively.

\\To check condition (\ref{teo}) we just need to  estimate $d_s$. Observe that $d_{1}= \D$ and
$$
d_{\b}\le \D {\D\choose\b}\le {\frac{\D^{1+\b}}{ \b!}}.
$$
Then, in the present case the function $\phi_{\FF}(\x)$ defined in (\ref{fea}) takes the form
$$
\phi_{\FF}(\x)\leq \frac{\D}{k} (\x+1) + \frac{1}{k^{\b}} {\frac{\D^{1+\b}}{ \b!}} (\x+1)^{\b}.
$$
%and we have that
%\be\label{kab}
%\inf_{x>0}\left[{\phi_{\EA}(x)\over x}\right] = \frac{1}{c^{\b}} {\D^{1+\b}\over \b!} \frac{\b^{\b}}{(\b-1)^{\b-1}}.
%\ee
%So, there is a $\b$-frugal coloring of $G$ as soon as
%\be\label{ex2aa}
%c\ge {\D^{1+ \frac{1}{\b}}\over {\b!}^{1/\b}} \b{(\b-1)^{\frac{1}{\b}-1}},
%\ee
And hence, with some calculation, we obtain the upper bound
$$
\ch_\b(G)\le {\D^{1+ \frac{1}{\b}}\over {\b!}^{1/\b}} \b{(\b-1)^{\frac{1}{\b}-1}} + \D,
$$
which, of course, is the same bound obtained in \cite{APS} via entropy compression method.

\end{document}